\documentclass[11pt]{article}

  \usepackage[letterpaper, left=1in, right=1in, top=1in, bottom=1in]{geometry}
  \providecommand{\keywords}[1]{\textbf{Keywords: } #1}
  
  \usepackage[affil-it]{authblk}  
  
  \usepackage{palatino}  
  \usepackage{setspace}  
  \usepackage{color}  
  \usepackage[dvipsnames]{xcolor}  

	\usepackage[labelfont=bf]{caption}	

  \usepackage{enumitem}  

  \usepackage{amsmath}  
  \usepackage{amssymb}  
  \usepackage{amsthm}  
  \usepackage{bm}  
  \usepackage{bbm}  
  \allowdisplaybreaks[4]  
	\usepackage{relsize}  
	\usepackage{algorithm}  
	\usepackage{algpseudocode}  
  
  \DeclareMathOperator*{\minimize}{minimize}
  
  \DeclareMathOperator*{\argmin}{arg\,min}
  \newcommand{\st}{\mathrm{subject\;to}}

	\theoremstyle{plain}

    \newtheorem{hypothesis}{Hypothesis}
	
  \usepackage{pifont}  

	\usepackage{longtable}  
	\usepackage{multirow}  
	\usepackage{array}  
	\usepackage{booktabs}  
    
  \usepackage{graphicx}  
  \usepackage{float}  
  \usepackage{subfig}  
\usepackage{tikz}
\usetikzlibrary{arrows.meta,positioning,shapes.geometric,calc}

	\usepackage{hyperref}
	\hypersetup{
		pdfstartview = {FitV},  
		colorlinks = true,    
		linkcolor = NavyBlue,  
		citecolor = NavyBlue,  
		filecolor = NavyBlue,  
		urlcolor = NavyBlue  
	}  

  \usepackage[round]{natbib}  

\usepackage{multicol}
\usepackage{caption}
\usepackage{subcaption}
\usepackage{enumitem,kantlipsum}
\usepackage{comment}
\usepackage{optidef}
\usepackage{alphalph}
\usepackage{etoolbox}
\usepackage{amsmath}
\usepackage{wrapfig}
\title{Uncovering expert objectives in production planning via inverse optimization: An industrial case study}
\author[1]{Shivi Dixit \thanks{These authors contributed equally to this work.}}
\author[1]{Rishabh Gupta \textsuperscript{*}}
\author[2]{Adam Kelloway}
\author[3]{John Wassick}
\author[1]{Qi Zhang \thanks{Corresponding author (qizh@umn.edu).}}
\affil[1]{Department of Chemical Engineering and Materials Science, University of Minnesota, Minneapolis, MN 55455, USA}
\affil[2]{The Dow Chemical Company, Midland, MI 48642, USA}
\affil[3]{Department of Chemical Engineering, Carnegie Mellon University, Pittsburgh, PA 15213, USA}
\date{}

\patchcmd{\subequations}{\alph{equation}}{\alphalph{\value{equation}}}{}{}

\begin{document}

\maketitle

\begin{abstract}
Production planning in the manufacturing industry often relies on the use of optimization models, but defining an appropriate objective function can be a challenge. In practice, planners must balance competing goals, manage uncertainty, and account for qualitative business preferences that are difficult to quantify. As a result, many optimization models fail to match expert behavior, limiting trust and adoption. In this work, we propose a data-driven inverse optimization framework to infer the objective function implicitly captured in expert planners’ decisions. We formulate the production planning problem as a mixed-integer linear program, where the unknown objective function is represented as a weighted sum of hypothesized cost terms. A suboptimality-loss-based inverse optimization method is then applied to learn the objective weights from historical production plans. The proposed approach is applied to a real industrial case provided by Dow, where the inferred weights reveal that avoiding inventory shortages and maintaining consistent cycle lengths dominate the planners’ decision-making. Time- and product-dependent extensions further improve predictive accuracy and uncover evolving priorities. Expert interviews confirm the practical validity of these insights. Overall, this study shows that inverse optimization can transform tacit human expertise into interpretable models, enabling more accurate and trusted decision-support tools for complex industrial systems.
\end{abstract}


\section{Introduction}

Production planning and scheduling is a critical function in chemical manufacturing and has been studied extensively in the literature \citep{harjunkoski2014scope, castro2018expanding}. Numerous optimization models have been proposed to tackle various production planning problems, typically formulated to minimize or maximize an objective function subject to relevant constraints. In practice, however, the appropriate objective function is often not explicitly known \citep{troutt2006behavioral}, which can be due to several reasons, such as: (i) There are hidden costs, such as those related to the impact of equipment wear and tear on plant reliability, that are difficult to quantify. (ii) Uncertainty, most notably in future product demand, affects planning decisions but is not modeled explicitly. (iii) It is unclear how competing objectives, such as cost, throughput, flexibility, and customer satisfaction, are prioritized. (iv) Planning decisions at a single plant can influence operations across the broader supply chain network, requiring non-intuitive adjustments to local objective functions to balance trade-offs at the network level. (v) The required information is simply not available due to a lack of documentation. When implementing production planning models in industrial applications, not knowing the right objective function can be a major challenge as it may lead to unreasonable solutions and ultimately the loss of user trust in the whole model.

Although the right objective function for a production planning problem may not be explicitly given, one can often assume that it is implicitly captured in the decisions made by an expert planner who has developed a deep intuition for the true costs and trade-offs through many years of experience. However, this information is not readily available in a form that can be directly incorporated into a mathematical model. Typically, it requires many discussions between the modeler and the planner and a considerable amount of trial and error to determine an objective function that provides good solutions. To facilitate this process, we develop in this work a systematic data-driven approach to learning the expert planner's mental objective function from production plans that the planner has generated in the past.

In the machine learning literature, the problem of learning how to perform a certain task by observing the behavior of a human expert is called \textit{imitation learning} \citep{hussein2017imitation}, where the goal is to find a mapping from inputs to actions such that these actions mimic the expert's behavior as closely as possible. It has found applications across a wide range of domains, most notably in robotics \citep{argall2009survey} and autonomous driving \citep{codevilla2018end}. \textit{Apprenticeship learning} \citep{abbeel2004apprenticeship} is a specific form of imitation learning that focuses on recovering or approximating the expert’s objective function rather than directly cloning their actions. As such, apprenticeship learning can result in more interpretable models as it directly helps to better understand how much the expert values certain aspects of the problem. It also tends to generalize better to unseen situations compared to other imitation learning methods \citep{zheng2022imitation}. 

The most commonly used method for apprenticeship learning is \textit{inverse reinforcement learning} (IRL). In standard reinforcement learning \citep{sutton1998reinforcement}, the decision-making problem is modeled as a Markov decision process with given possible states and actions, the corresponding transition probabilities, and a known reward function; the goal is then to learn a policy, which maps states to actions, that maximizes the cumulative reward. IRL addresses the inverse problem where the reward function is unknown but can be inferred from observed state–action trajectories, assuming that they result from maximizing an underlying cumulative reward. It is especially well suited for problems where it is difficult to model the underlying process, e.g. a helicopter flying under various weather conditions, in detail. IRL has received significant attention in the past two decades in which many methodological advancements have been made, focusing on improving model accuracy, data efficiency, and generalizability. For a recent review of IRL, we refer the reader to \citet{arora2021survey}.

A method closely related to IRL is \textit{inverse optimization} (IO) \citep{Ahuja2001}. Here, the decision-making problem is modeled as a mathematical optimization problem, where typically the constraints are known but the objective function is not. Similar to IRL, the goal of IO is to learn the underlying objective function (or reward function) given observed expert decisions that are assumed to be optimal or near-optimal solutions to the optimization problem. Early IO studies focused on \textit{classical formulations} that enforced a perfect fit for a single observation, most often in linear or network flow models \citep{Burton1992, Ahuja2001, Heuberger2004}. In contrast, more recent \textit{data-driven IO} frameworks accommodate multiple and noisy observations by minimizing measures of suboptimality across data \citep{Keshavarz2011, MohajerinEsfahani2018a, Aswani2018}. This evolution has made IO particularly useful for learning the implicit objectives of expert decision-makers, whose real-world choices often reflect heuristics or approximate optimality. Methodological advances have further expanded IO to nonlinear \citep{gupta2022efficient, lu2025bo4io}, discrete \citep{Wang2009, bulut2021complexity, bodur2022inverse}, and high-dimensional models, with scalable algorithms now available to estimate objective functions from observed behavior \citep{tan2020learning, gupta2021, gupta2022efficient, bodur2022inverse}. For a comprehensive overview of the IO research landscape, we refer the reader to \citet{ChanReviewPaper}.

In recent years, IO has been applied across diverse domains. In healthcare, it has been used to infer the relative weight clinicians place on competing treatment goals in radiation therapy planning, sometimes integrated with machine learning for greater interpretability \citep{babier2021ensemble, ajayi2022objective}. In energy systems, IO has helped estimate consumer demand response preferences and utility curves from smart grid data \citep{saez2017short, esteban2024estimating}. Transportation studies have employed IO to recover congestion cost functions and route-choice preferences from observed traffic flows \citep{zhang2017data, zattoni2025inverse, chan2024conformal}. In finance, IO has been applied to infer investors’ risk-aversion coefficients from portfolio holdings \citep{alsabah2021robo, yu2023learning}. More recently, it has also been used to uncover institutional objectives in financial aid allocation \citep{bara2025revealing} and to analyze the implicit priorities behind political gerrymandering \citep{smith2024gap}.

Several of these studies have leveraged real-world data: for example, clinical treatment plans in radiation oncology \citep{ajayi2022objective}, consumption data from electricity markets \citep{saez2017short}, traffic flow observations \citep{zhang2017data}, and actual investment portfolios \citep{alsabah2021robo}. These empirical validations demonstrate IO’s ability to extract meaningful objective functions that explain expert or system-level decision-making. However, a common feature across these applications is that the forward problems typically assume relatively simple constraints such as linear dose-volume limits in healthcare, flow conservation in transportation, or budget constraints in finance. While this simplification facilitates tractability and interpretability, it underutilizes IO’s potential. As we show in this work, incorporating richer, domain-specific constraints allows IO models to capture more nuanced decision processes, thereby increasing both their explanatory power and practical relevance.

In this work, we apply data-driven IO to an industrial case study provided by The Dow Chemical Company (subsequently referred to as Dow), where we uncover expert planners' objectives when performing production planning for a particular plant. Here, the underlying planning problem can be readily formulated as a mixed-integer linear program (MILP) for which we do not initially have an appropriate objective function. In this case, with the help of an expert planner, we can identify potential objectives that may play a role in the planners' decision-making; however, we do not know to what extent they drive planning decisions and how the planners balance the trade-offs between competing objectives. We assume that the true objective function is a weighted sum of multiple objective terms and use IO to infer the weights from historical production plans generated by the planners. As a result, the trained model is inherently interpretable as the learned weights directly reveal the relative importance of the various hypothesized objectives. After obtaining the inferred objective function, we conduct an interview with an expert planner to evaluate the results and gain additional insights into how well the model aligns with the planners' true decision-making process. To the best of our knowledge, this is the first study reported in the literature that applies data-driven IO to real-world production planning data in an effort to elucidate and quantify human expert planners' objectives. It effectively demonstrates the efficacy of the IO approach in complex industrial applications.

The remainder of this paper is organized as follows. In Section \ref{sec:Problem}, we provide details about the production planning problem, the corresponding MILP formulation, and the specific form of the objective function that we are trying to learn. Then, the proposed data-driven IO approach is described in Section \ref{sec:Method}. In Section \ref{sec:Results}, we present the results from our real-world industrial case study and discuss these results based on feedback and insights received from an expert planner. Finally, we close with concluding remarks in Section \ref{sec:Conclusions}.

\section{Problem description}
\label{sec:Problem}

In this work, we consider a specific production planning problem from Dow; however, the setup is representative of many planning applications across the process industry such that the proposed approach can be readily extended to other similar problems. In the following, we provide a statement of the production planning problem, a mathematical formulation of the scheduling constraints, and a list of hypothesized objectives that the planners may have been considering in their decision-making.

\subsection{Production planning problem}

\begin{figure}[ht!]
    \centering
    \includegraphics[width=0.75\linewidth]{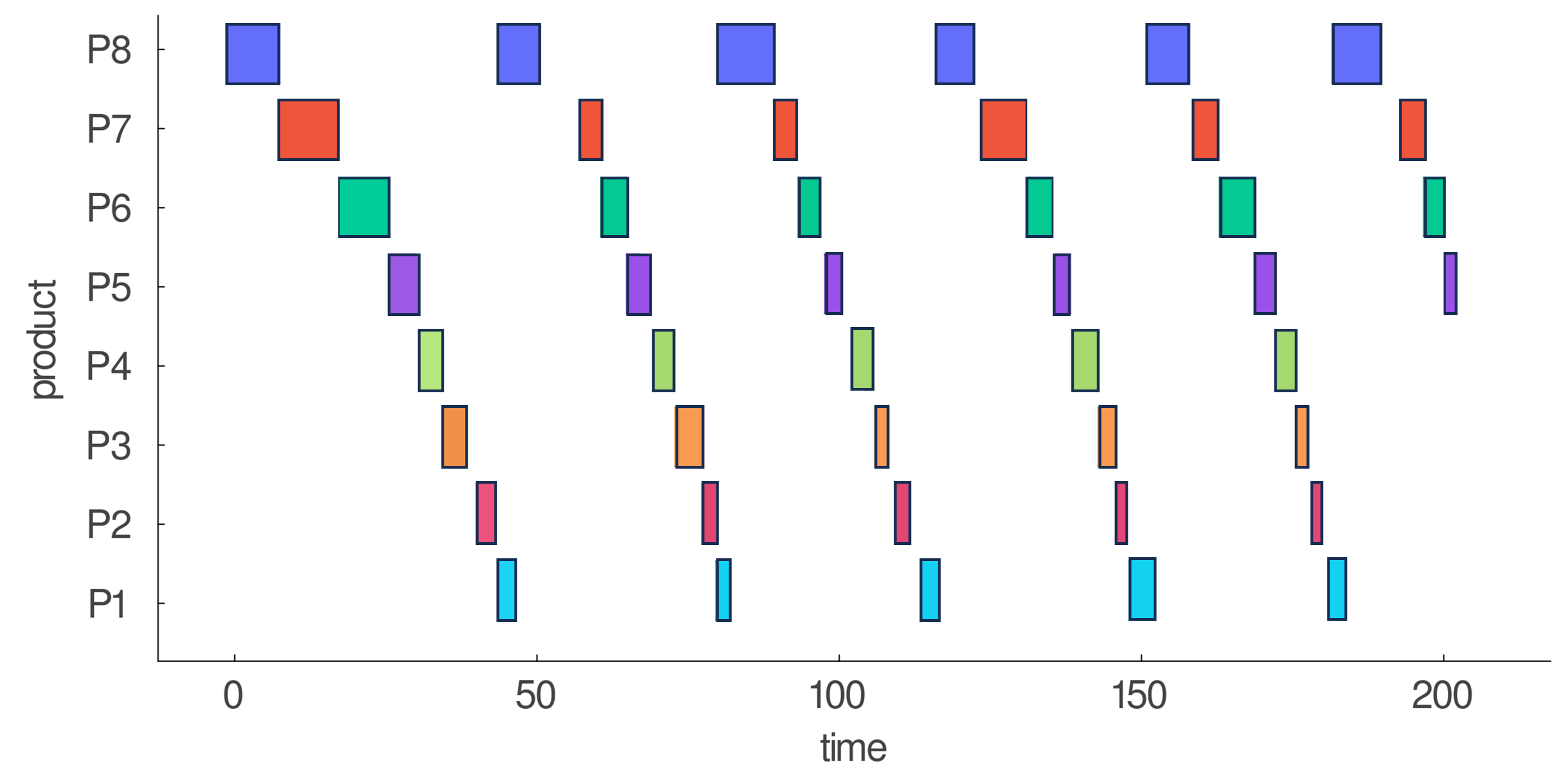}
    \caption{Example of a Gantt chart depicting a typical production plan, where each bar indicates a campaign for a particular product.}
    \label{fig:Gantt}
\end{figure}

In this problem, we consider a plant that manufactures multiple products to meet demand over a given planning horizon, which generally requires multiple production cycles, where each cycle consists of multiple production campaigns, one for each product. As illustrated in the example Gantt chart shown in Figure \ref{fig:Gantt}, where each bar depicts a campaign, the order in which the products are produced is the same in all cycles. This production sequence is determined in a higher-level planning problem that considers a longer time horizon and other aspects such as transition costs and off-grade products. With the sequencing of campaigns being fixed, this problem focuses on optimizing the lengths of the campaigns, which can vary across different production cycles. Here, a campaign consists of a set of consecutive batches of the same product; hence, the campaign length depends on the number of batches produced in that campaign. Campaigns cannot overlap since there is only one production unit, and there is a minimum number of batches for each campaign so that no product is skipped in a cycle. Moreover, there are restrictions on the length of each production cycle and on the total number of cycles in the planning horizon.

\begin{figure}[ht!]
    \centering
    \includegraphics[width=0.65\linewidth]{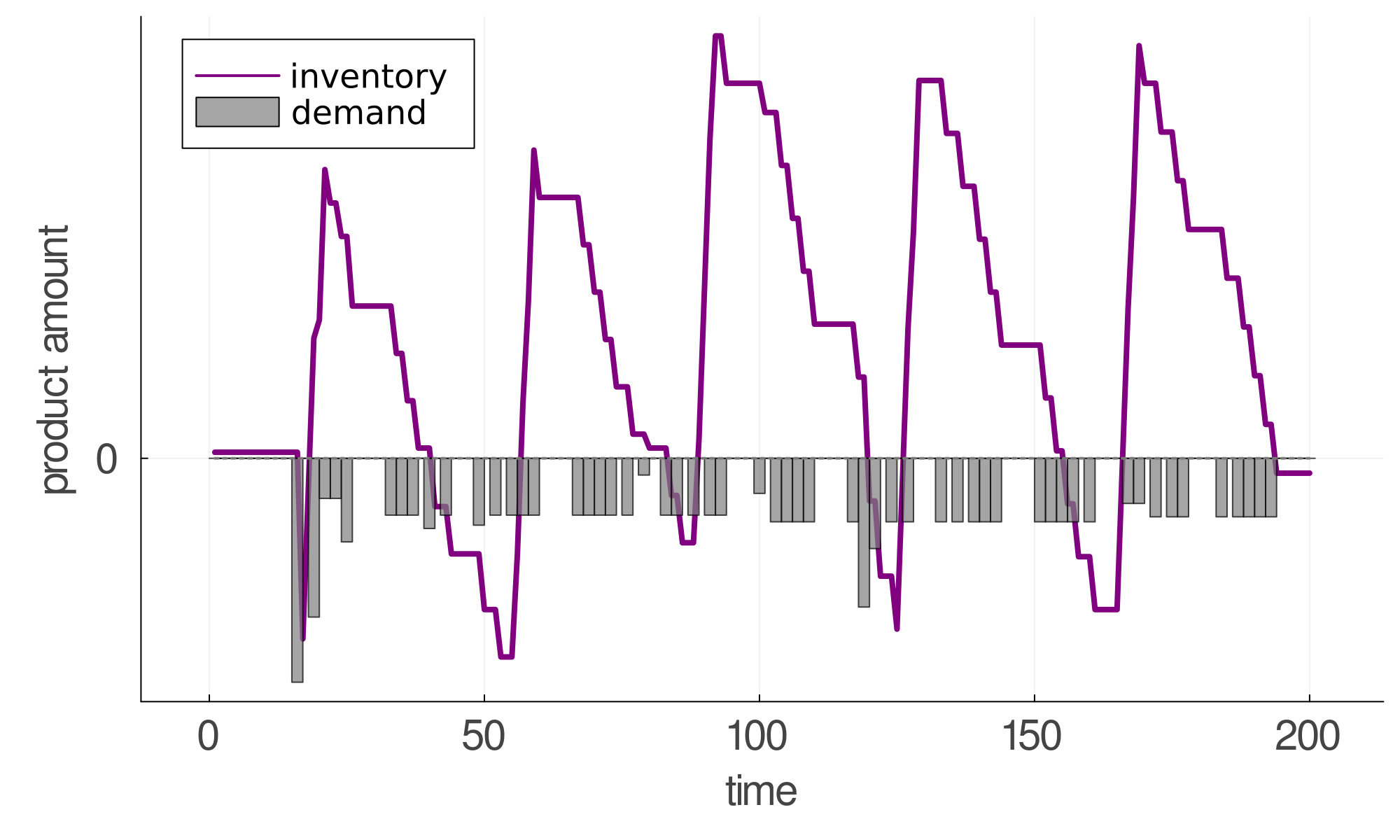}
    \caption{Example of typical demand and inventory profiles for a particular product.}
    \label{fig:Inventory}
\end{figure}

Figure \ref{fig:Inventory} shows a typical demand profile and the corresponding inventory profile for a particular product. Note that here demand values are given as negative numbers, as the product is drawn from the storage tank to meet demand. The inventory profile clearly reflects the timing of the campaigns. During a campaign, the inventory level of the corresponding product typically increases, while it decreases afterward as the inventory is depleted to meet demand until the next campaign for the same product starts. In this case, inventory levels can be negative, which indicate backlogs, i.e. demands that are not met in the current cycle but can be satisfied in later cycles. Lower and upper bounds are generally imposed on the inventory levels.

The given demand profile for the planning problem considers both firm customer orders and forecast demand. As orders are usually only known for the next few time periods, the production plan needs to be updated on a regular basis. This means that we only implement the part of the production decisions that apply to the time periods before the next time the production plan is re-optimized. It is important to keep this rolling-horizon planning setup in mind as it may have an impact on how human planners make decisions.

\subsection{Scheduling constraints}

We now describe the mixed-integer constraints that jointly define the set of feasible solutions for the production planning problem. We apply a discrete-time formulation with a planning horizon of $T$ time periods, where the notation is such that time period $t$ starts at time $t-1$ and ends at time $t$. The resulting inventory balance constraints are as follows:
\begin{align}
    & I_{pt} = I_{p,t-1} + \beta \, x_{pt} - D_{pt} \quad \forall \, p \in \mathcal{P}, \, t \in \mathcal{T} \label{eqn:inv_balance} \\
    & I^{\min}_p \leq I_{pt} \leq I^{\max}_p \quad \forall \, p \in \mathcal{P}, \, t \in \mathcal{T} \label{eqn:inv_bounds} \\
    & x_{pt} \in \{0,1\} \quad \forall \, p \in \mathcal{P}, \, t \in \mathcal{T},
\end{align}
where $\mathcal{P} = \{1,\dots,P\}$ denotes the set of $P$ products ordered according to the given production sequence, $\mathcal{T}=\{1,\dots,T\}$ denotes the set of time periods, $D_{pt}$ is the demand for product $p$ in time period $t$, and $\beta$ is the fixed batch size, which is the amount of product manufactured in one time period. Note that in our particular problem, the batch size is the same for all products; however, if needed, $\beta$ can be product-dependent. The continuous variable $I_{pt}$ denotes the inventory level of product $p$ at time $t$ and can take negative values to account for backlogs in meeting customer demands, and the binary variable $x_{pt}$ equals 1 if product $p$ is produced in time period $t$ and 0 otherwise. Equations \eqref{eqn:inv_balance} simply state that the inventory level at time $t$ is the inventory level at time $t-1$ plus the amount produced in time period $t$ minus the amount withdrawn to meet demand in time period $t$. Constraints \eqref{eqn:inv_bounds} enforce minimum and maximum inventory levels denoted by $I^{\min}_p$ and $I^{\max}_p$, respectively.

In this problem setting, only one production unit is available. This means that only one product can be produced at a time, which is captured in the following constraints:
\begin{equation}
    \sum_{p \in \mathcal{P}} x_{pt} \leq 1 \quad \forall \, t\in \mathcal{T} .
\end{equation}

The next set of constraints are used to track the timing and length of each campaign. We introduce the binary variables $y_{pt}$ and $z_{pt}$ that indicate the start and end, respectively, of a campaign producing product $p$, which is captured in the following constraints:
\begin{align}
   & 2 \, y_{p,t-1}-1 \leq x_{pt}-x_{p,t-1} \quad \forall \, p \in \mathcal{P}, \, t \in \mathcal{T} \label{eqn:start1} \\
   & y_{p,t-1} \geq x_{pt}-x_{p,t-1} \quad \forall \, p \in \mathcal{P}, \, t \in \mathcal{T} \label{eqn:start2} \\
   & 2 \, z_{p,t-1}-1 \leq x_{p,t-1}-x_{pt} \quad \forall \, p \in \mathcal{P}, \, t \in \mathcal{T} \label{eqn:end1} \\
   & z_{p,t-1} \geq x_{p,t-1}-x_{pt} \quad \forall \, p \in \mathcal{P}, \, t \in \mathcal{T} \label{eqn:end2} \\
   & y_{pt} \in \{0,1\}, \; z_{pt} \in \{0,1\} \quad \forall \, p \in \mathcal{P}, \, t = \{0,\dots,T\}.
\end{align}
Here, constraints \eqref{eqn:start1} and \eqref{eqn:start2} ensure that $y_{pt}$ equals 1 if a campaign producing product $p$ starts at time $t$ and 0 otherwise. Similarly, constraints \eqref{eqn:end1} and \eqref{eqn:end2} ensure that $z_{pt}$ equals 1 if a campaign producing product $p$ ends at time $t$ and 0 otherwise.

The length of a campaign producing product $p$ must be between $\gamma^{\min}_p$ and $\gamma^{\max}_p$ time periods, where we assume that $1 \leq \gamma^{\min}_p \leq \gamma^{max}_p$. The following constraints enforce these bounds:
\begin{align}
    & 1-y_{pt} \geq \sum_{t'=t}^{\min \{t+\gamma^{\min}_p-1, T\}} z_{pt'} \quad \forall \, p \in \mathcal{P}, \, t \in \{0,.., T-1\} \\
    & y_{pt} \leq \sum_{t'=t+1}^{t+\gamma^{\max}_p} z_{pt'} \quad \forall \, p \in \mathcal{P}, \, t \in \{0,.., T-\gamma^{\max}_p \}.
\end{align}

The following constraints restrict the total number of cycles within the planning horizon and the length of each cycle:
\begin{align}
    & \zeta^{\min} \leq \sum_{t \in \mathcal{T}} y_{pt} \leq \zeta^{\max} \quad \forall \, p \in \mathcal{P} \label{eqn:cycle_number} \\
    & 1-y_{pt} \geq \sum_{t'=t+1}^{\min\{t+\eta^{\min}-1,T\}} y_{pt'} \quad \forall \, p \in \mathcal{P}, \, t \in \{0,.., T-1\} \text{ when } \eta^{\min} \geq 2 \label{eqn:cycle_len_min} \\
    & y_{pt} \leq \sum_{t'=t+\eta^{\min}}^{t+\eta^{\max}} y_{pt'} \quad \forall \; p \in \mathcal{P}, \; t \in \{0,..,T-\eta^{\max}\}, \label{eqn:cycle_len_max}
\end{align}
where constraints \eqref{eqn:cycle_number} ensure that the number of cycles is between $\zeta^{\min}$ and $\zeta^{\max}$. Note that the last cycle can be incomplete in the sense that there may not be enough time remaining in the planning horizon to cycle through all products in that cycle. Constraints \eqref{eqn:cycle_len_min} and \eqref{eqn:cycle_len_max} state that the length of each cycle must be between $\eta^{\min}$ and $\eta^{\max}$ time periods, where we have $\sum\limits_{p \in \mathcal{P}} \gamma^{\min}_p \leq \eta^{\min} \leq \eta^{\max}$.

Finally, we need the following constraints to ensure that the given production sequence is followed in each production cycle:
\begin{align}
    & z_{pt} \leq \sum_{t'=t}^{t+\theta_p} y_{p+1,t'} \quad \forall \, p \in \{1,\dots,P-1\}, \, t \in \{1,.., T-\theta_p \} \label{eqn:sequence1} \\
    & z_{P,t} \leq \sum_{t'=t}^{t+\theta_P} y_{1,t'} \quad \forall \, t \in \{1,.., T-\theta_P \}, \label{eqn:sequence2}
\end{align}
where $\theta_p$ denotes the maximum number of time periods between the end of a campaign producing product $p$ and the start of the next campaign. We assume that $\theta_p < \eta^{\min} - \gamma^{\max}_p$ and $\theta_p < \sum\limits_{p' \in \mathcal{P} \setminus \{p\}} \gamma^{\min}_{p'} + P - 1$. The first condition ensures that the next campaign for the next product starts before the next campaign for the same product, which will take place in the next cycle. The second condition ensures that we cannot have a production gap after every campaign. In practice, these conditions are not overly restrictive as $\theta_p$ can still be reasonably large so that it is not limiting; we usually observe production gaps, if any, that are much shorter in duration.

Note that throughout the proposed model formulation, we do not make use of a separate index for cycles or campaigns since each campaign is uniquely determined by its start (as indicated by $y_{pt}$) and end (as indicated by $z_{pt}$). By doing so, we avoid having to prespecify the number of cycles and vary it in multiple runs of the model to find the optimal number of cycles.

\subsection{Hypothesized objectives}
\label{sec:Objectives}

In this problem, we assume that the scheduling constraints are known, which is reasonable as they represent physical constraints and well-established business rules for which most parameter values are readily available, but we do not know the objective function that aligns with expert planners' decisions. To infer such an objective function, we first interviewed a planner from Dow to gain an understanding of the main factors that drive their decision-making. It turns out that it is not that difficult to determine a list of objectives that planners may care about; however, many of these objectives are competing with each other, and it is unclear to what extent the planners emphasize one over another.

In the following, we present a list of hypothesized objectives derived from our discussion with the planner and show how they can be modeled. We then propose a single objective function that combines these objectives and identify the parameters that need to be estimated. Note that there are some additional business-specific objectives that we could have incorporated; however, for confidentiality reasons, we only consider the more generally applicable objectives in our case study. As we will show in our results, just considering these potential objectives was enough to provide useful insights with our proposed approach.

\paragraph{Minimizing inventory holding cost.} The one objective that immediately comes to mind is minimizing inventory holding cost, a cost component that is incorporated in most production planning models. The corresponding objective function term can be written as follows:
\begin{equation}
    C_{ic} = \sum_{p \in \mathcal{P}} \sum_{t \in \mathcal{T}} \phi_p \max\{0, I_{pt}\},
\end{equation}
where $\phi_p$ is the unit inventory holding cost for product $p$, and $C_{ic}$ is the total inventory cost, where $ic$ is a short hand notation for inventory cost.

\begin{figure}[ht!]
    \centering
    \includegraphics[width=0.6\linewidth]{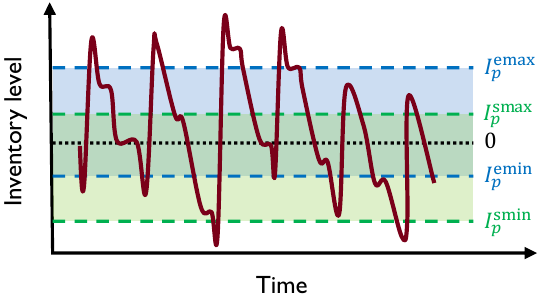}
    \caption{Illustration of desired inventory ranges at the start and end of each campaign.}
    \label{fig:inv_bounds}
\end{figure}

\paragraph{Keeping inventory levels within given bounds.} Production planners typically aim to maintain sufficient on-hand inventory to meet demand while avoiding excessive stock levels that lead to high inventory holding costs. In this case, the planners want to keep the inventory levels within some desired ranges, which are illustrated in Figure \ref{fig:inv_bounds}. At the start of a campaign, the inventory level for a product $p$ should preferably be between $I^{\mathrm{smin}}_p$ and $I^{\mathrm{smax}}_p$, while at the end of a campaign, it should be between $I^{\mathrm{emin}}_p$ and $I^{\mathrm{emax}}_p$. 

Given these inventory targets, it makes sense to penalize inventory levels that are outside these ranges. We can calculate these deviations using the following constraints:
\begin{align}
    & M (y_{pt}-1) + I^{\mathrm{smin}}_p - I_{pt} \leq S^{l}_{pt} \quad \forall \, p \in \mathcal{P}, \, t \in \mathcal{T} \\
    & M (y_{pt}-1) + I_{pt} - I^{\mathrm{smax}}_p \leq S^{u}_{pt} \quad \forall \, p \in \mathcal{P}, \, t \in \mathcal{T} \\
    & M (z_{pt}-1) + I^{\mathrm{emin}}_p - I_{pt} \leq E^{l}_{pt} \quad \forall \, p \in \mathcal{P}, \, t \in \mathcal{T} \\
    & M (z_{pt}-1) + I_{pt} - I^{\mathrm{emax}}_p \leq E^{u}_{pt} \quad \forall \, p \in \mathcal{P}, \, t \in \mathcal{T} \\
    & S^l_{pt} \geq 0, \; S^u_{pt} \geq 0, \; E^l_{pt} \geq 0, \; E^u_{pt} \geq 0 \quad \forall \, p \in \mathcal{P}, \, t \in \mathcal{T},
\end{align}
where $M$ is a sufficiently large positive parameter. The variables $S^l_{pt}$ and $S^u_{pt}$ are used to compute the positive deviations from the lower and upper bounds $I^{\mathrm{smin}}_p$ and $I^{\mathrm{smax}}_p$, respectively, at the start of a campaign. Similarly, $E^l_{pt}$ and $E^u_{pt}$ represent the deviations at the end of a campaign. For each type of inventory level deviation, we obtain an additional objective function term as follows:
\begin{align}
    & C_{sl} = \sum_{p \in \mathcal{P}} \sum_{t \in \mathcal{T}} S^l_{pt} \\
    & C_{su} = \sum_{p \in \mathcal{P}} \sum_{t \in \mathcal{T}} S^u_{pt} \\
    & C_{el} = \sum_{p \in \mathcal{P}} \sum_{t \in \mathcal{T}} E^l_{pt} \\
    & C_{eu} = \sum_{p \in \mathcal{P}} \sum_{t \in \mathcal{T}} E^u_{pt}.
\end{align}

\paragraph{Maintaining consistent cycle lengths.} Planners also try to schedule the production such that the lengths of the production cycles are close to a desired cycle length previously determined at a higher level of the planning decision hierarchy. Specifically, they want the length of each cycle to be within a given range $[\delta^{\min}, \delta^{\max}]$, where $\eta^{\min} \leq \delta^{\min} \leq \delta^{\max} \leq \eta^{\max}$. We can use the following constraints to identify cycles whose lengths are outside of this range:
\begin{align}
    & y_{pt} \leq \sum_{t'=t+\delta^{\min}}^{t+\delta^{\max}} y_{pt'} + v_t \quad \forall \, p \in \mathcal{P}, \, t \in \{0,..,T-\delta^{\max}\} \\
    & v_t \in \{0,1\} \quad \forall \, t \in \{0,..,T-\delta^{\max}\},
\end{align}
where the binary variable $v_t$ must equal 1 if the number of time periods between the start of a campaign that begins at time $t$ and the start of the next campaign producing the same product is less than $\delta^{\min}$ or greater than $\delta^{\max}$. Summing $v_t$ over all applicable time periods provides an objective function term that can be used to encourage consistency in cycle times:
\begin{equation}
    C_{cl} = \sum_{t=0}^{T-\delta^{\max}} v_t.
\end{equation}

\paragraph{Minimizing plant downtime.} Apart from planned shutdowns, planners generally prefer to run the plant continuously with as few production gaps as possible. To model this objective, we introduce the following constraints:
\begin{align}
    & z_{pt} \leq y_{p+1,t} + w_{pt} \quad \forall \, p \in \{1,\dots,P-1\}, \, t \in \{1,.., T-1 \} \\
    & z_{P,t} \leq y_{1,t} + w_{P,t} \quad \forall \, t \in \{1,.., T-1 \} \\
    & w_{pt} \in \{0,1\} \quad \forall \, p \in \mathcal{P}, \, t \in \{1,.., T-1 \},
\end{align}
where the binary variable $w_{pt}$ must equal 1 if there is a gap in production after a campaign producing product $p$ ends at time $t$. The objective is then to minimize the total number of production gaps:
\begin{equation}
    C_g = \sum_{p \in \mathcal{P}} \sum_{t=1}^{T-1} w_{pt} \label{eqn:C_g}.
\end{equation}

\paragraph{Weighted-sum objective function.} Given the hypothesized objectives defined above, we posit the planners' true objective function to be a weighted sum of these objective function terms, i.e. we assume that the underlying optimization problem can be formulated as the following MILP:
\begin{equation}
\tag{P}
\label{eqn:optproblem}
\begin{aligned}
    \minimize \quad & \sum_{k \in \mathcal{K}} \alpha_k \, \rho_k \, C_k \\
    \st \quad & \text{constraints \eqref{eqn:inv_balance}--\eqref{eqn:C_g}},
\end{aligned}    
\end{equation}
where $\mathcal{K} = \{ic, sl, su, el, eu, cl, g\}$ is the set of hypothesized objectives, $\rho_k$ is a scaling factor that is calculated beforehand and ensures the framework is not biased toward specific terms when the order of magnitudes vary significantly between different cost terms, and $\alpha_k$ is the weighting factor for objective $k$. One can interpret this as the planners solving a multi-objective optimization problem, where the weighting factors can directly tell us how the planners balance the different objectives to obtain trade-off solutions that they deem optimal. However, we do not know the values of the weighting factors; hence, our goal is to infer them from production plans previously generated by the planners. Note that although $\alpha$ are unknown, they are parameters in problem \eqref{eqn:optproblem}, which is an MILP.

\section{Data-driven inverse optimization approach}
\label{sec:Method}

We pose the problem of uncovering expert planners’ unknown preferences, i.e. inferring the weighting factors $\alpha$ in \eqref{eqn:optproblem}, as a data-driven IO problem. Specifically, since the underlying optimization model is an MILP, it becomes an inverse mixed-integer linear optimization problem (IMILOP). IMILOPs have been studied from a theoretical and algorithmic standpoint \citep{Wang2009, bodur2022inverse}, and more recently through data-driven formulations that accommodate multiple and noisy observations \citep{moghaddass2021inverse, holani2026}. The latter are particularly relevant for our setting, as human planning decisions often deviate from exact optimality due to uncertainty and operational heuristics. 

Conceptually, our approach in this work is most closely related to the IMILOP framework proposed by \citet{moghaddass2021inverse}, who introduced a suboptimality-minimization formulation to handle imperfect decisions. 
Here, we extend this idea to enable the recovery of an interpretable objective function directly from real production planning data. In the following, we describe the forward and inverse formulations, the suboptimality-based reformulation used to handle noisy data, and the cutting-plane solution method employed to solve the resulting semi-infinite program efficiently.

\subsection{Forward optimization problem}

As detailed in the previous section, we model the expert planners' decision process as the MILP given in \eqref{eqn:optproblem}, which we now write in the following compact form for notational convenience:
\begin{equation}
\tag{FOP}
\label{eqn:FOP}
\begin{aligned}  
\minimize_{x \in \mathbb{R}^p \times \mathbb{Z}^q} \quad & c^{\top} x \\
\st \quad & A(u)x \leq b(u), \\
\end{aligned}
\end{equation}
where $x$ are the mixed-integer decision variables, $u$ are the input parameters (e.g. customer demands and current inventory levels), and $(A(u), b(u))$ capture the linear constraints, which change with $u$. The cost vector $c$ defines the linear objective function; as such, it includes the weighting factors from formulation \eqref{eqn:optproblem} that we are trying to learn.

\subsection{Inverse optimization problem}

IO aims to infer $c$ from observed input–decision pairs. Let $\mathcal{I}$ denote the set of historical instances that serve as training data points, with $(u_i, x_i)$ denoting the inputs and the corresponding observed decisions for instance $i \in \mathcal{I}$. The IMILOP can then be formulated as follows:
\begin{equation}
\tag{IOP}
\label{eqn:IOP}
\begin{aligned}  
\underset{\hat{c} \in C, \, \hat{x}}{\minimize} \quad & \sum_{i \in \mathcal{I}} \|x_i - \hat{x}_i\| + \lambda \|\hat{c}\|_1 \\
\st \quad & \hat{x}_i \in \argmin_{\tilde{x} \in \mathbb{R}^p \times \mathbb{Z}^q} \{\, \hat{c}^{\top}\tilde{x} : A(u_i)\tilde{x} \leq b(u_i) \} \quad \forall \, i \in \mathcal{I},
\end{aligned}
\end{equation}
where $\hat{x}_i$ denotes the predicted decisions using $\hat{c}$ for instance $i$. As indicated by the objective function of \eqref{eqn:IOP}, the goal is to minimize the difference between the observed and predicted decisions and a regularization term that encourages the construction of sparse models. The constraints state that $\hat{x}_i$ is an optimal solution to the corresponding MILP. Also, we use $\mathcal{C}$ to denote the set from which $\hat{c}$ can be chosen, which may include constraints that, for example, enforce nonnegativity and prohibit the trivial solution where all cost coefficients are zero.

\subsection{Reformulation via suboptimality minimization}

Problem \eqref{eqn:IOP} is a bilevel optimization problem with $|\mathcal{I}|$ lower-level problems. Because these lower-level problems are MILPs, we cannot apply KKT- or duality-based single-level reformulations. We could apply a value-function reformulation and a cutting-plane method to solve the resulting semi-infinite program; however, one key challenge remains, namely that we still need to enforce that every predicted solution $\hat{x}_i$ satisfies all mixed-integer constraints, which significantly increases the computational complexity of the problem as the number of training data points increases. 

To address the computational challenge, we apply the following alternative formulation of the IMILOP proposed by \citet{moghaddass2021inverse}, which minimizes a suboptimality loss instead of the decision loss shown in \eqref{eqn:IOP}:
\begin{equation}
\tag{IOP-SL}
\label{eqn:IOP-SL}
\begin{aligned}  
\minimize_{\hat{c} \in C, \, \epsilon \in \mathbb{R}_+^{|\mathcal{I}|}} \quad & \sum_{i \in \mathcal{I}} \epsilon_i + \lambda \|\hat{c}\|_1 \\
\st \quad & \hat{c}^{\top} x_i - \hat{c}^{\top}\tilde{x}_i \leq \epsilon_i \quad \forall \, i \in \mathcal{I}, \,\tilde{x}_i \in \mathcal{S}_i,
\end{aligned}
\end{equation}
where $\mathcal{S}_i := \{x \in \mathbb{R}^p \times \mathbb{Z}^q : A(u_i)x \leq b(u_i)\}$ denotes the feasible region of instance $i$. Here, we introduce the nonnegative variable $\epsilon_i$ which must be positive if $x_i$ is not the optimal solution that minimizes the cost defined by $\hat{c}$ in instance $i$. As problem \eqref{eqn:IOP-SL} minimizes the sum of all $\epsilon_i$'s, they provide a measure of the level of suboptimality for all observed decisions given a chosen $\hat{c}$. Note that no mixed-integer constraints are involved in \eqref{eqn:IOP-SL} since it does not include any predicted solutions that must satisfy the constraints of the forward problem.

\subsection{Solution method}

\begin{figure}[ht!]
\centering
\begin{tikzpicture}[
  node distance=10mm and 16mm,
  >=Latex,
  every node/.style={font=\small},
  box/.style={rectangle, rounded corners, draw, align=center, minimum width=34mm, minimum height=8mm},
  decision/.style={diamond, draw, aspect=2, inner sep=1.2pt, align=center},
  io/.style={trapezium, trapezium left angle=70, trapezium right angle=110, draw, align=center, minimum width=36mm, minimum height=8mm}
]

\node[box] (start) {Start};
\node[box, below=of start] (init) {Initialize $\bar{\mathcal{S}}_i$ for all $i \in \mathcal{I}$};
\node[io, below=of init] (master) {Solve master problem};
\node[box, below=of master] (cand) {Estimated cost vector $\hat{c}^*$};
\node[io, below=of cand, text width=55mm] (subs) {Solve cut-generating problem with $\hat{c}$ for each $i \in \mathcal{I}$ to obtain $x_i^\star$};
\node[box, below=of subs] (update) {Set $\bar{\mathcal{S}}_i \leftarrow \bar{\mathcal{S}}_i \cup \{x_i^\star\}$ for all $i \in \mathcal{I}$};
\node[decision, below=of update] (check) {Converged?};
\node[box, below=of check] (done) {Output $\hat{c}^*$ as the final estimate};

\draw[->] (start) -- (init);
\draw[->] (init) -- (master);
\draw[->] (master) -- (cand);
\draw[->] (cand) -- (subs);
\draw[->] (subs) -- (update);
\draw[->] (update) -- (check);
\draw[->] (check.east) -- ++(25mm,0) |- (master.east);
\draw[->] (check) -- node[right, xshift=1mm]{Yes} (done);
\node[above right=5mm and 26mm of check.east] (lbl) {No};
\end{tikzpicture}
\caption{Cutting-plane scheme for solving problem \eqref{eqn:IOP-SL}: alternate between solving the master problem on a finite set of cuts and generating new cuts via the cut-generating problems.}
\label{fig:io-cutting-plane-flow}
\end{figure}
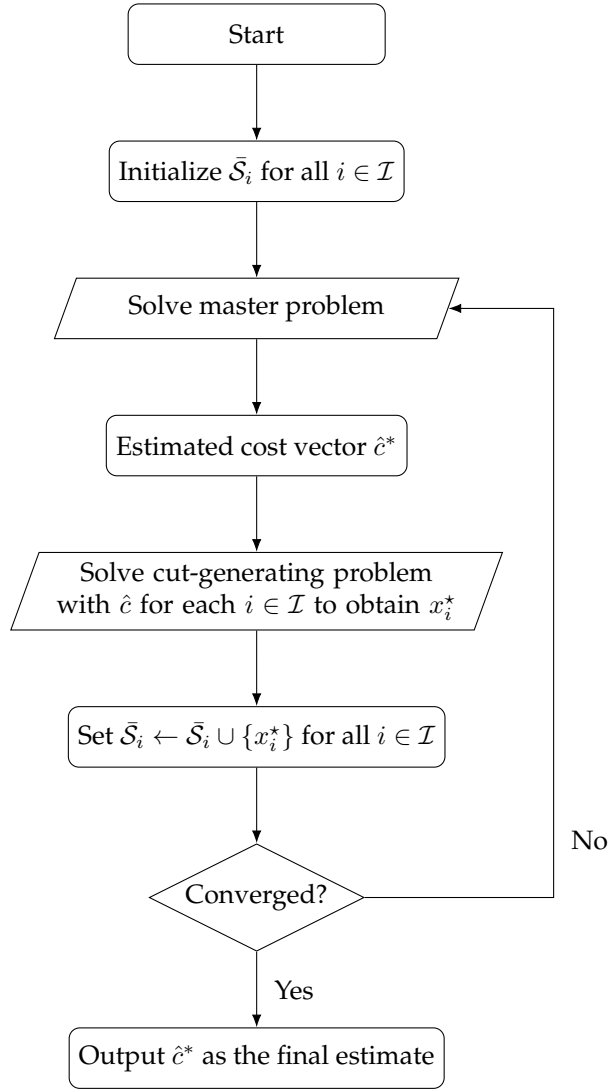

Problem \eqref{eqn:IOP-SL} is a semi-infinite program since the set $\mathcal{S}_i$ generally contains an infinite number of points. We solve it using a cutting-plane approach that alternates between a master problem and a set of cut-generating subproblems.
\begin{itemize}
    \item The \textit{master problem} is a relaxation of \eqref{eqn:IOP-SL} where each $\mathcal{S}_i$ is replaced by a finite subset $\bar{\mathcal{S}}_i \subseteq \mathcal{S}_i$, which results in a finite number of constraints. As a result, the master problem is just a linear program (LP). Solving it provides a candidate cost vector $\hat{c}^*$ and $\epsilon^*_i$ for each $i \in \mathcal{I}$.
    \item Given a $\hat{c}^*$ obtained from solving the master problem, we can formulate a \textit{cut-generating problem} for each $i \in \mathcal{I}$ that is \eqref{eqn:FOP} with $c \gets \hat{c}^*$ and $u \gets u_i$. By solving this problem, we aim to find solutions that potentially achieve cost values lower than the one for the observed decisions with the current estimates of the cost coefficients $\hat{c}^*$ minus $\epsilon^*_i$; for any $i \in \mathcal{I}$ where this applies, we add the corresponding $x^*_i$ to $\bar{\mathcal{S}}_i$.
\end{itemize}
The cutting-plane algorithm, as illustrated in the flowchart shown in Figure \ref{fig:io-cutting-plane-flow}, terminates when $(\hat{c}^*)^{\top} (x_i - x^*_i) \leq \epsilon^*_i$ for all $i \in \mathcal{I}$ or the change in the objective function value of \eqref{eqn:IOP-SL} or in the value of $\hat{c}^*$ is small over multiple consecutive iterations or a maximum number of iterations is reached. Note that the cut-generating problem only needs to be solved to optimality in the final iterations to show convergence.

\section{Results, discussion, and qualitative validation}
\label{sec:Results}

We now apply the proposed framework to a real industrial case provided by Dow. In the discussion of our results, we emphasize the inherent interpretability of the estimated cost parameters, which allowed us to derive specific hypotheses about the planners' decision-making process. Using these interpretable results as a basis, we sought feedback from an expert planner, who was able to qualitatively confirm or invalidate our hypotheses. In the following, we detail our analysis and discussion. We list all our hypotheses, and for each one, we assess its validity and provide an explanation that resulted from our conversation with the expert planner. Note that due to confidentiality reasons, we cannot disclose all plant and product specifications as well as the actual product demands. Therefore, all results presented in the following are given without units, and the values are normalized if necessary. However, this anonymization does not hinder the analysis and interpretation of the results.

\subsection{Data and model training}

The plant considered in this case study produces eight products. There are multiple planners responsible for scheduling the production at this plant. We assume that their decisions are driven by the same underlying objective function; hence, we treat production plans generated by different planners as if they originated from the same decision-maker. We are provided with 70 historical production plans, each for a planning horizon of 200 time periods. It is important to note that each plan represents the planner's decisions at a given point in time on the \textit{planned} actions over the given time horizon; it is not what eventually got \textit{implemented} (except for the decisions in the first few time periods) since the production plan is re-optimized on a regular basis in a rolling-horizon fashion. 

Of the 70 available production plans, we use 50 for training and 20 for testing. Models are constructed and evaluated using five different randomly chosen training-testing splits to assess the robustness of the proposed approach and quantify the variability in the parameter estimates and model predictions. Hence, when means and variances are reported in the results, they originate from these five different training-testing splits. Problem \eqref{eqn:IOP-SL} is solved using the cutting-plane algorithm with a maximum number of iterations of 20. We adaptively increase the maximum solution time for the cut-generating problem with each iteration to allow quick cut generation in the earlier iterations and determine convergence toward the end as the cut-generating problem is solved to optimality in the final iterations. 
Once a model is trained, i.e. estimates of the unknown cost coefficients are obtained from IO, we can predict production planning decisions by solving problem \eqref{eqn:optproblem}.

All optimization problems were implemented in Julia v1.7.2 using the modeling environment JuMP v0.22.3 \citep{DunningHuchetteLubin2017} and solved with  Gurobi v10.3. All instances were solved utilizing 24 cores and 60 GB of memory on the Mesabi cluster of the Minnesota Supercomputing Institute (MSI) equipped with Intel Haswell E5-2680v3 processors.


\subsection{Model selection}

We first compare the performance of three different models that differ in the number of cost terms considered in the objective function:
\begin{itemize}
    \item a one-parameter model (model-1) that only considers inventory holding costs, i.e. all weighting factors in the objective function of problem \eqref{eqn:optproblem} are set to zero except for $\alpha_{ic}$;
    \item a five-parameter model (model-5) that, in addition to inventory holding costs, also includes the cost terms related to the four different types of penalties on inventory levels being outside of the desired ranges, i.e. also $\alpha_{sl}$, $\alpha_{su}$, $\alpha_{el}$, and $\alpha_{eu}$ can be nonzero;
    \item and the full seven-parameter model (model-7) that incorporates all seven objective function terms described in Section \ref{sec:Objectives}, including the penalties on cycle length deviations and production gaps with $\alpha_{cl}$ and $\alpha_g$ being their respective weighting factors.
\end{itemize}

\begin{figure}[ht!]
    \centering
    \includegraphics[width=10cm]{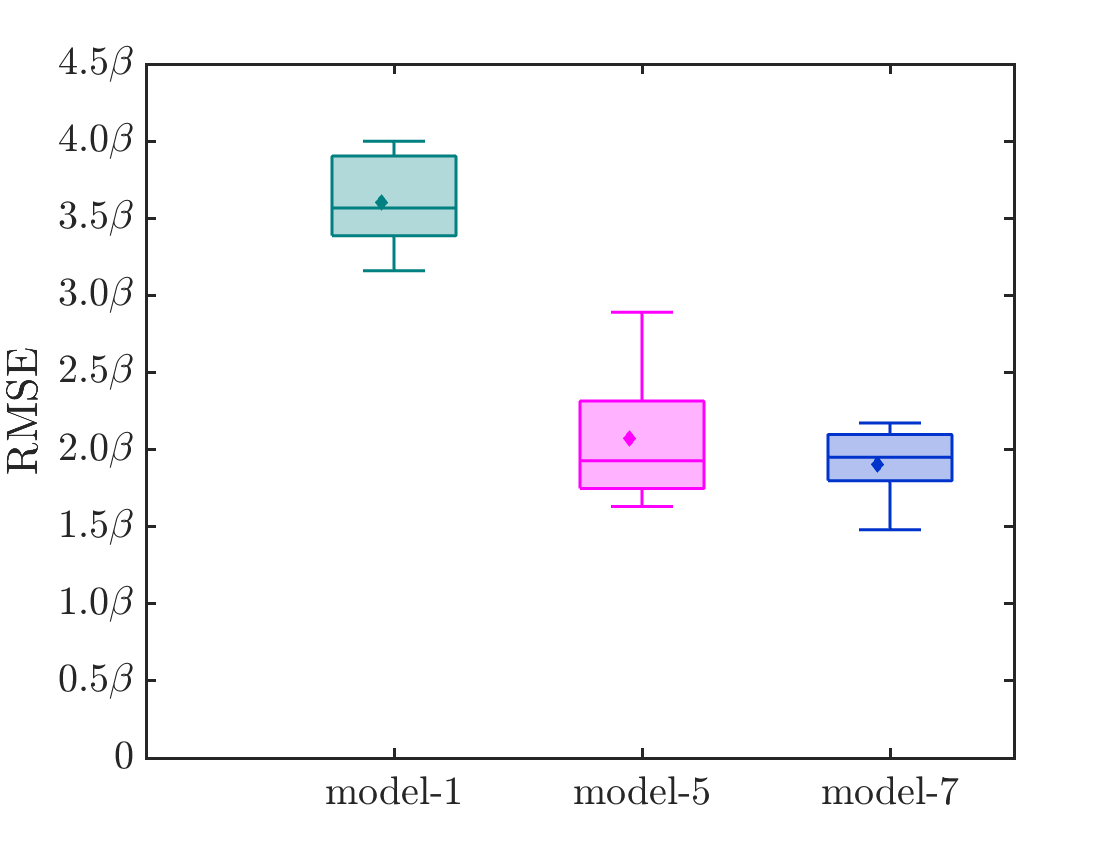}
    \caption{RMSEs achieved by model-1, model-5, and model-7.}
    \label{fig:rmse_overall}
\end{figure}

To obtain model-1, we simply set $\alpha_{ic}$ to 1 since there is no relative difference in importance between multiple objectives to be learned. We use model-1 to provide a baseline as its objective function captures the only directly known financial costs, i.e. the inventory holding costs. As such, it is the model one would initially formulate in the absence of any information about other costs that are not directly quantifiable, such as those capturing the other six hypothesized objectives. Model-5 and model-7 are trained using the proposed data-driven IO approach; we use them to assess the impact of assuming different sets of objective function terms. Figure \ref{fig:rmse_overall} compares the root mean squared errors (RMSEs) achieved by these three models on the test dataset, where one can see that the prediction accuracy increases with the number of model parameters. Note that the RMSEs are computed based on the difference between the actual and predicted inventory profiles and are given in terms of the batch size $\beta$.

\begin{figure}[hb!]
    \centering
    \subfloat{\includegraphics[width=5.0cm]{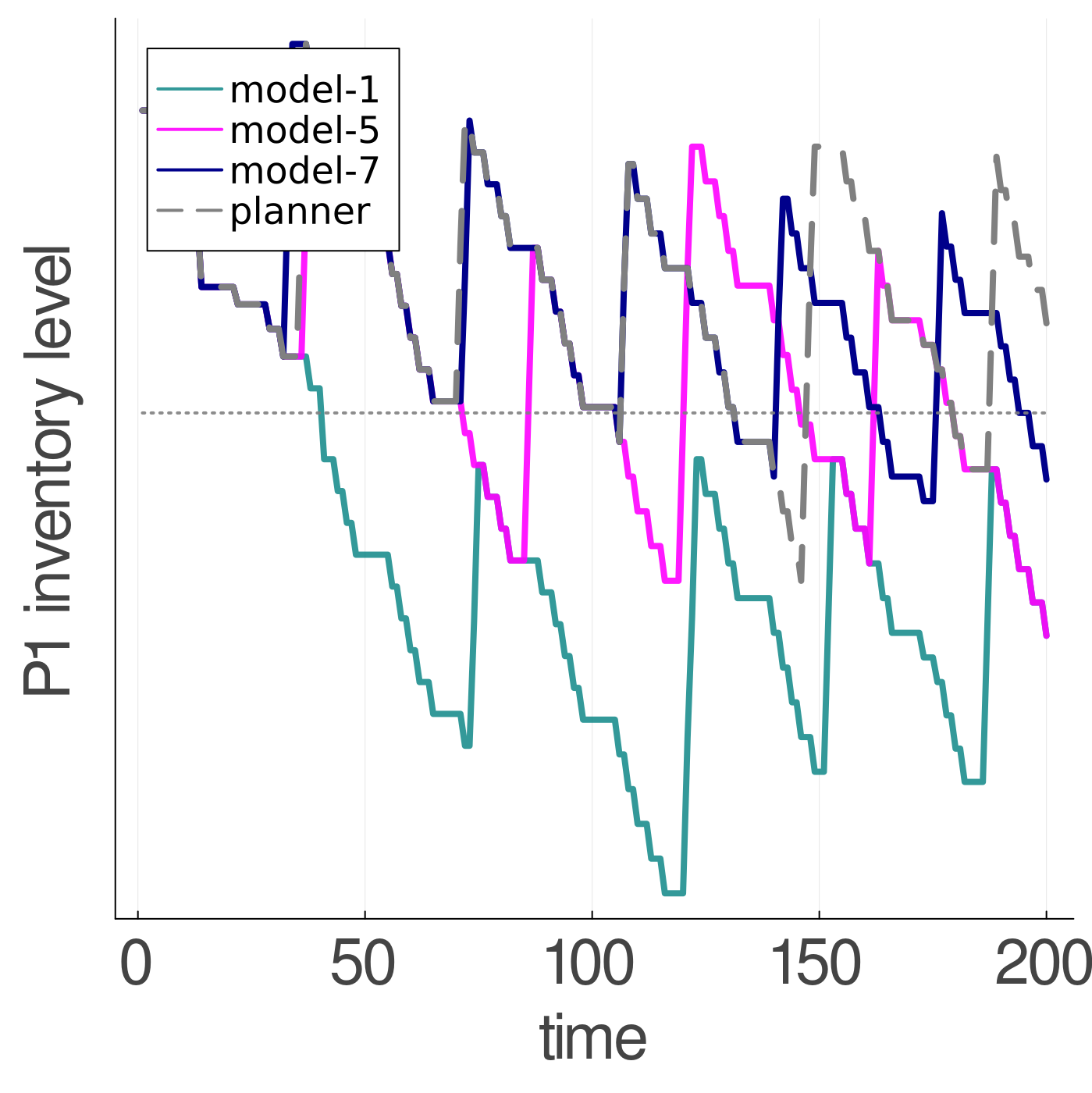}}
    \subfloat{\includegraphics[width=5.0cm]{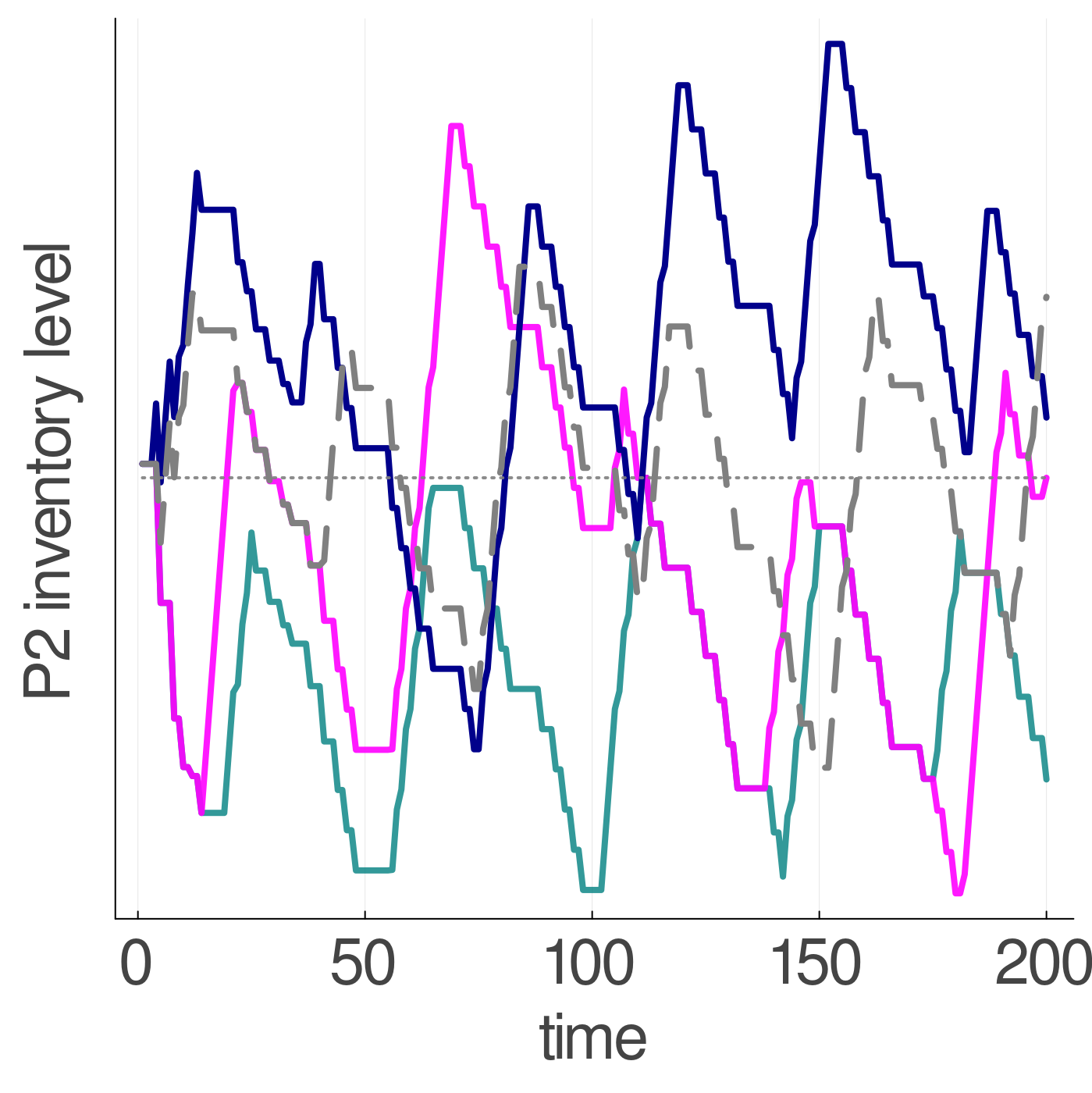}}
    \subfloat{\includegraphics[width=5.0cm]{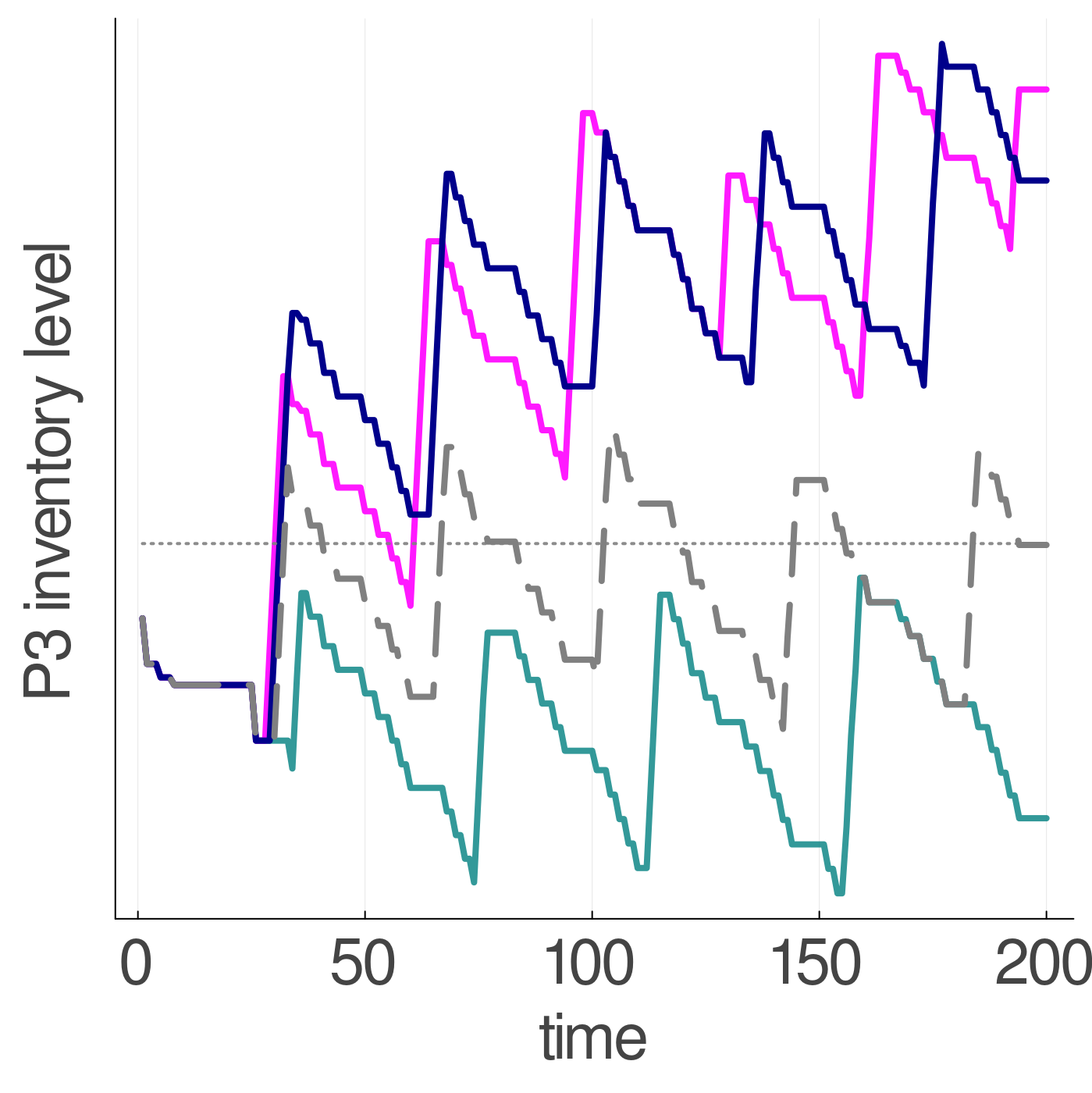}} \\
    \subfloat{\includegraphics[width=5.0cm]{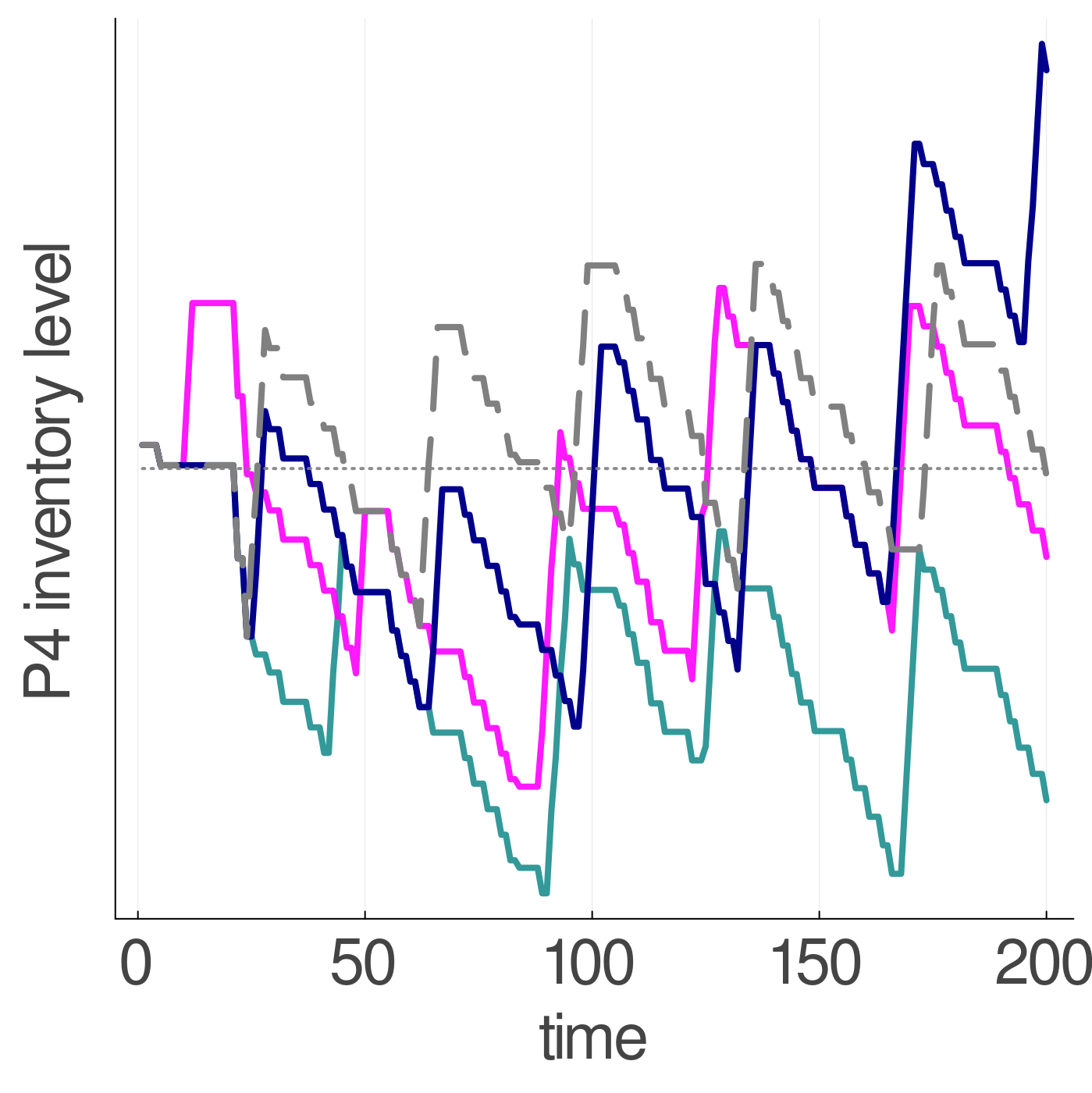}}
    \subfloat{\includegraphics[width=5.0cm]{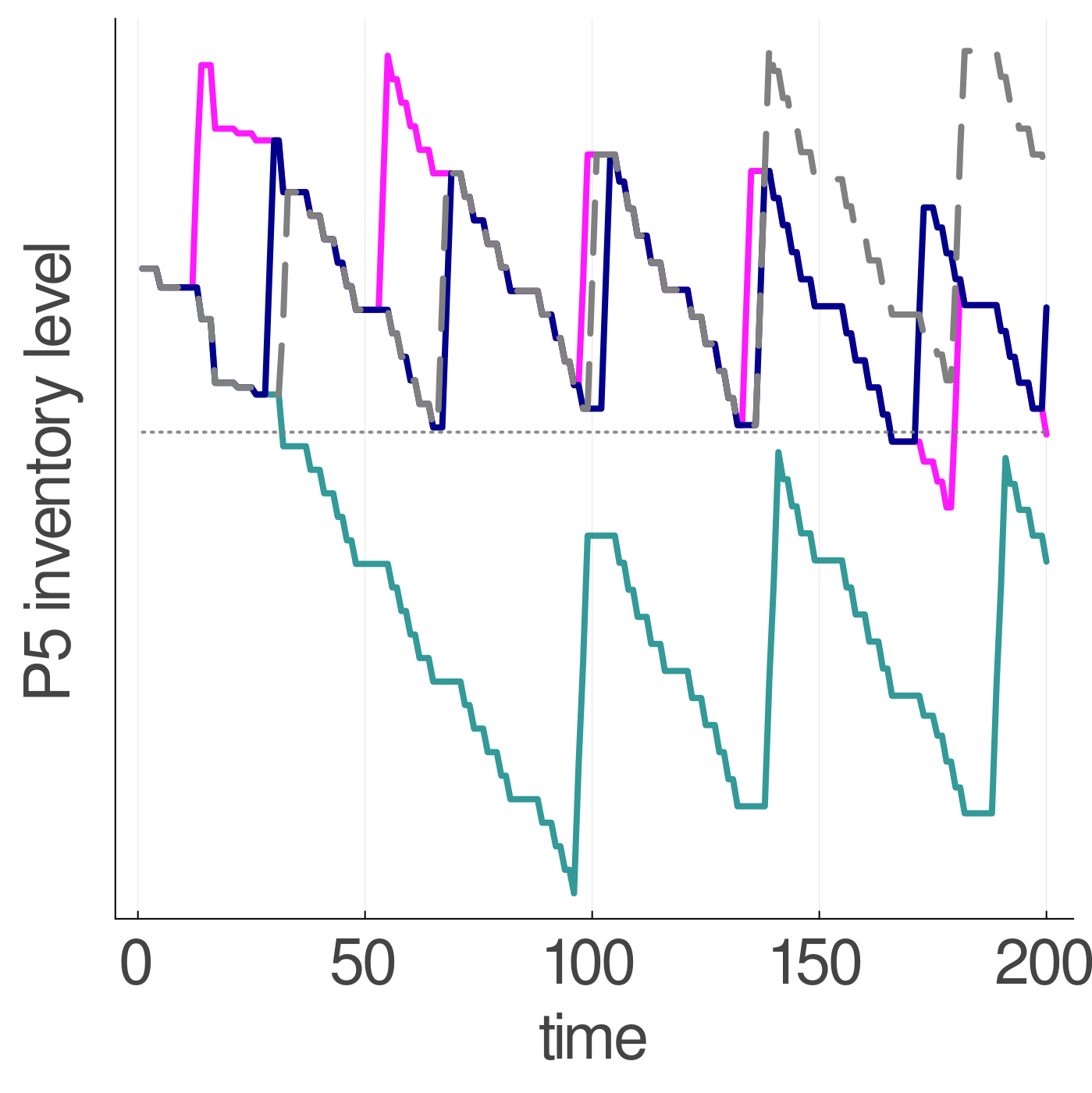}}
    \subfloat{\includegraphics[width=5.0cm]{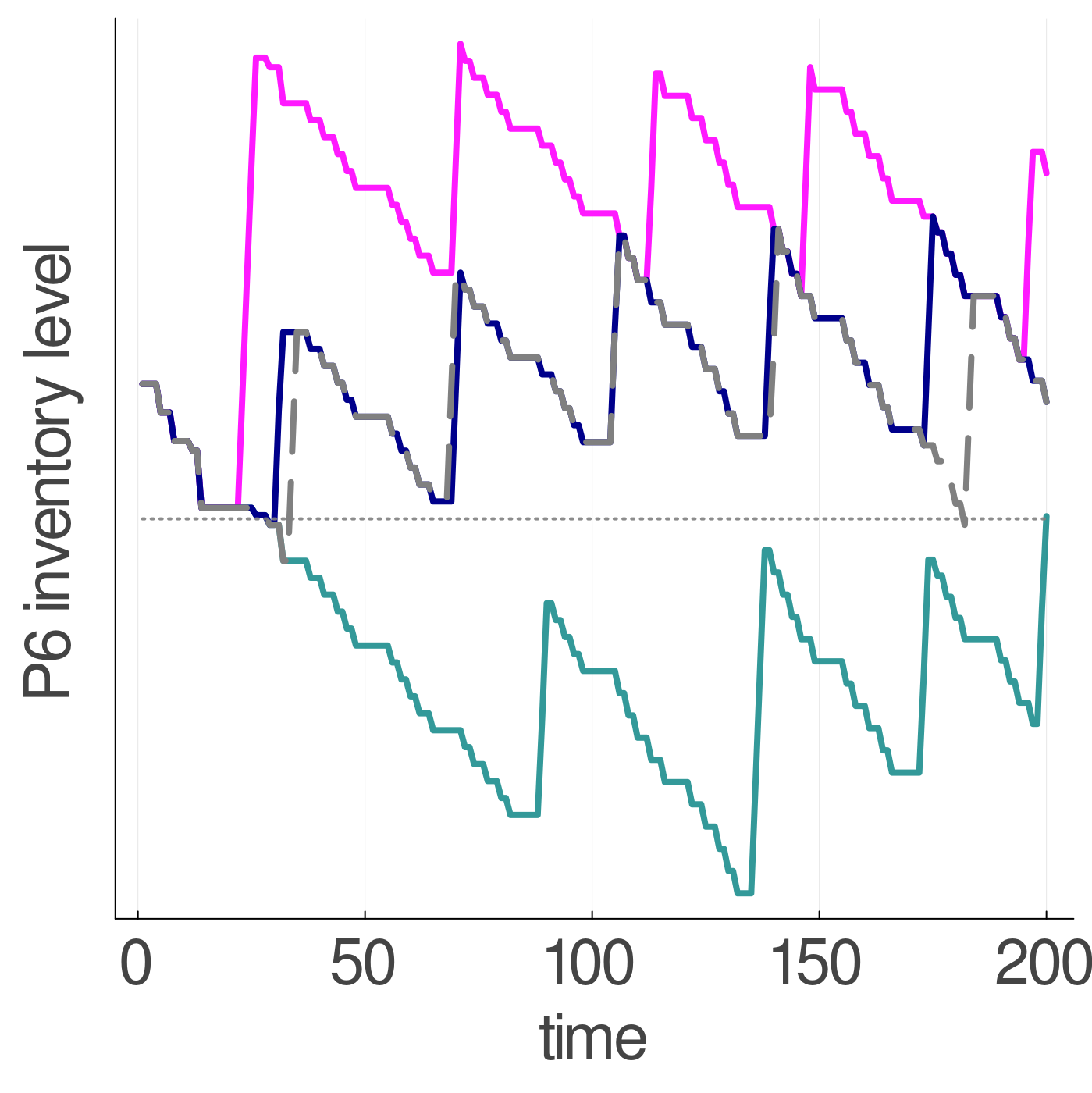}} \\
    \subfloat{\includegraphics[width=5.0cm]{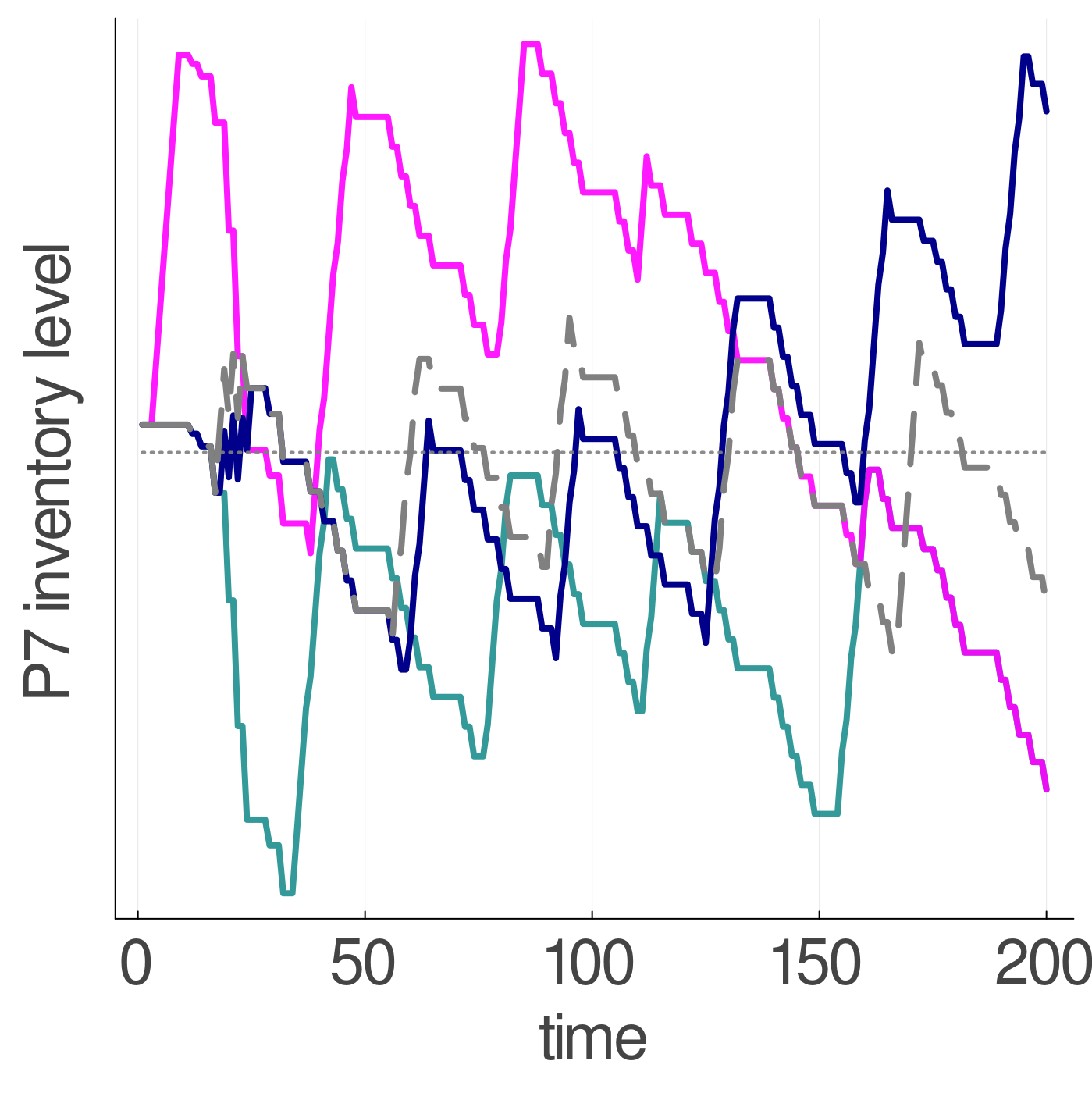}}
    \subfloat{\includegraphics[width=5.0cm]{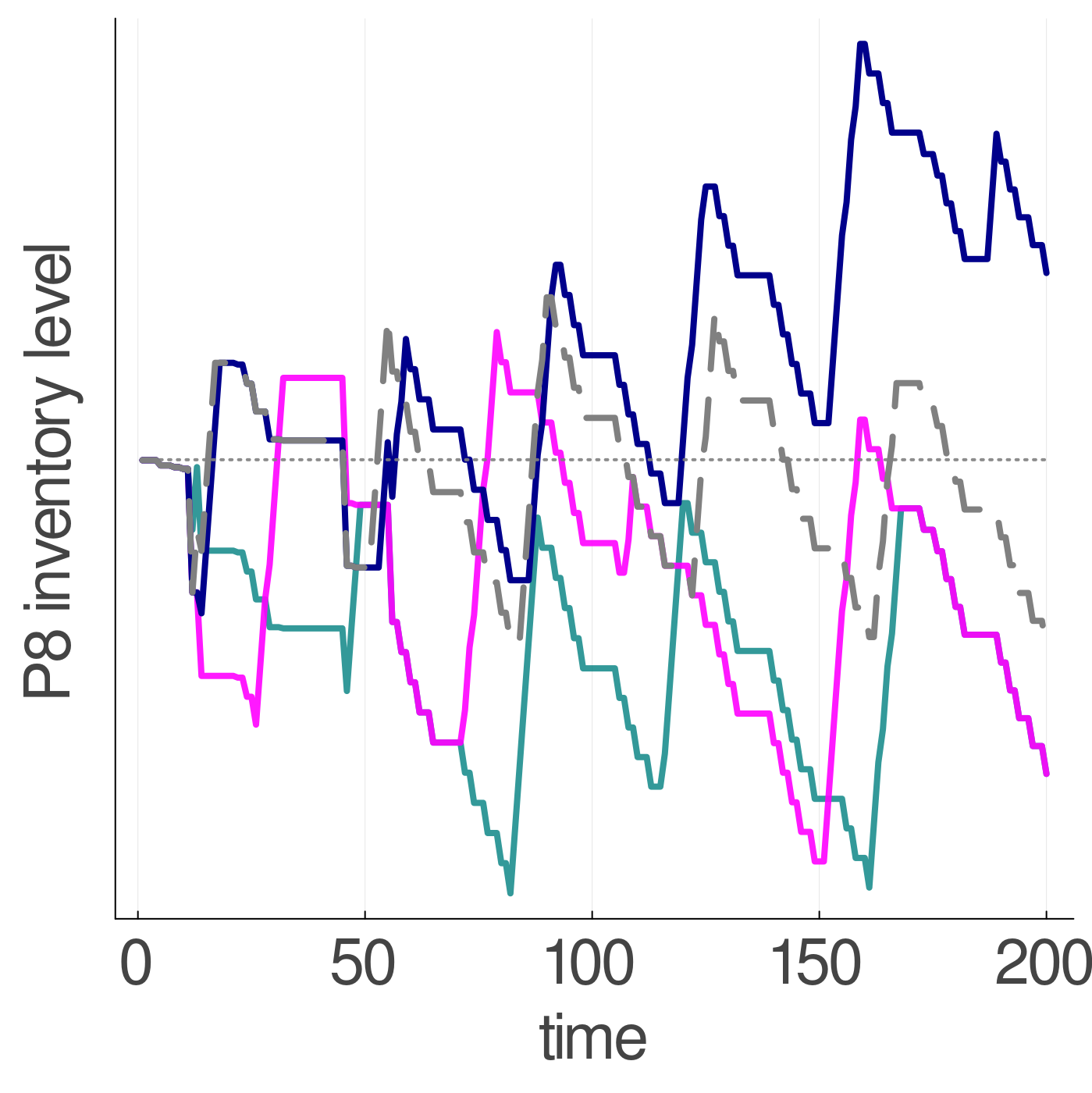}}
    \caption{Comparing for a particular instance the actual inventory profiles for all eight products with the profiles predicted by model-1, model-5, and model-7.}
    \label{fig:T200}
\end{figure}

Using a particular instance, we examine the differences in the predicted production schedules, which are shown in terms of the resulting inventory profiles in Figure \ref{fig:T200}. We can see that compared to the planners' actual decisions, model-1 consistently predicts much lower inventory levels while model-5 and model-7 achieve more realistic inventory profiles. This indicates that inventory holding costs, which are the only costs in this problem that can be readily quantified, are not the sole driver of the planners' decisions.

\subsection{Analysis and interpretation}

\begin{figure}[ht!]
    \centering
    \includegraphics[width=15cm]{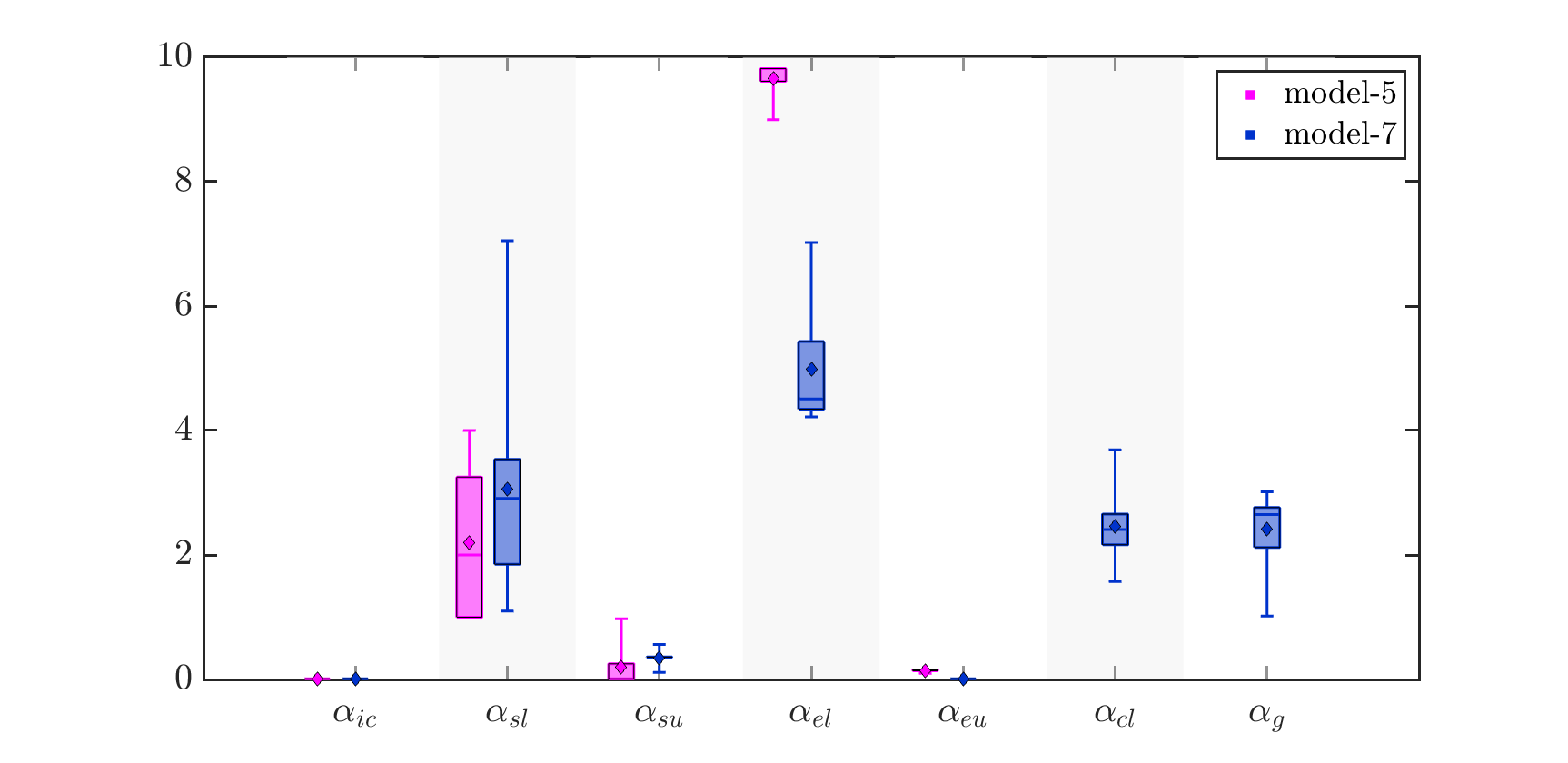}
    \caption{Estimates of the weighting factors obtained for model-5 and model-7.}
    \label{fig:costpref}
\end{figure}

We now examine the estimates of the weighting factors $\alpha$, which are shown in Figure \ref{fig:costpref} for model-5 and model-7. Recall that the objective function terms in problem \eqref{eqn:optproblem} are normalized using the scaling factors $\rho$; this allows us to gain insights into the relative importance between different objectives by directly comparing the corresponding $\alpha$-values. In Figure \ref{fig:costpref}, one can see that while the values of the first five weighting factors are different in the two models, the ranking of these parameters according to their values is the same for both models. This indicates that both models capture the same relative importance across these five objectives. However, since model-7 with the two additional parameters provides more accurate predictions, we use the model-7 results as a basis for the following discussion.

The first interesting and initially surprising observation we make is that $\alpha_{ic} = 0$, which leads to the following hypothesis.

\begin{hypothesis}
    Planners do not directly consider inventory holdings costs to inform their planning decisions. \textup{-- True. According to the expert planner we spoke to, the planners never explicitly compute the inventory holding costs in their decision-making process; instead, they manage inventory by focusing on other guidelines such as keeping inventory levels within desired ranges. This does not mean that inventory costs are not important, but rather that planners do not solely focus on reducing inventory since they also need to balance this with the benefits of having sufficient on-hand inventory to maintain a high service level.}
\end{hypothesis}

Comparing the values of $\alpha_{sl}$, $\alpha_{su}$, $\alpha_{el}$, and $\alpha_{eu}$, we see that $\alpha_{sl} \gg \alpha_{su}$ and $\alpha_{el} \gg \alpha_{eu}$. This indicates that the planners penalize inventory levels below the given ranges more heavily than those above.

\begin{hypothesis}
    It is more important to avoid understock than overstock. \textup{-- True. Planners typically would rather have excess inventory than not enough as they do not want to miss out on revenue opportunities when demand is higher than expected; the downside of having higher inventory holding costs when demand turns out to be lower is less severe.}
\end{hypothesis}

Next, we observe that $\alpha_{sl} < \alpha_{el}$, which we interpret as follows.

\begin{hypothesis}
\label{hyp:el_more_important}
    Planners pay more attention to the inventory level at the end rather than at the start of a campaign. \textup{-- False. Planners look at both the start and the end of a campaign. If there were a preference of one over the other, it would be the opposite; the rationale is that one should pay more attention to the start of a campaign since there is more opportunity to affect the end of a campaign as the production plan gets updated over time.}
\end{hypothesis}

The values of $\alpha_{cl}$ and $\alpha_g$ are relatively large, which leads to the following two hypotheses.

\begin{hypothesis}
    It is important to maintain consistent cycle lengths. \textup{-- True. Doing so strikes a balance between avoiding overly long production cycles to keep inventory cost low and having sufficiently long cycles to hedge against demand uncertainty.}
\end{hypothesis}

\begin{hypothesis}
    There is a strong desire to reduce plant downtime. \textup{-- True. Planners typically try to run the plant continuously without any idle times unless there is a planned shutdown for maintenance purposes. }
\end{hypothesis}

We see that $\alpha_g$ has a value similar to $\alpha_{cl}$ and is less than $\alpha_{sl}$ and $\alpha_{el}$. From this observation, we further derive the following more nuanced hypothesis.

\begin{hypothesis}
    Planners need to carefully balance the benefits of reducing production gaps and keeping inventory levels within the desired bounds. \textup{-- False. According to the expert planner, it is very rare that there is a significant benefit from idling the plant so that it is usually not a consideration in the planners' decision-making. In other words, under normal circumstances, planners do not intentionally create gaps in their production plans even at the risk of high inventory costs.}
\end{hypothesis}

Given the insight provided by the expert planner, one would expect $\alpha_g$ to be significantly larger. A potential reason for this not being the case is that we do actually observe a few, mostly quite short, production gaps in the provided production plans, which we are trying to match in the our predictions; hence, the model cannot place a very high penalty that would effectively disallow any gaps in production. However, we do not know exactly why these plant downtimes occurred; maybe these were due to some external factors that are not included in our model but also not considered by the planners under normal circumstances.

\subsection{Time- and product-dependent objectives}

\begin{figure}[ht!]
    \centering
    \includegraphics[width=15cm]{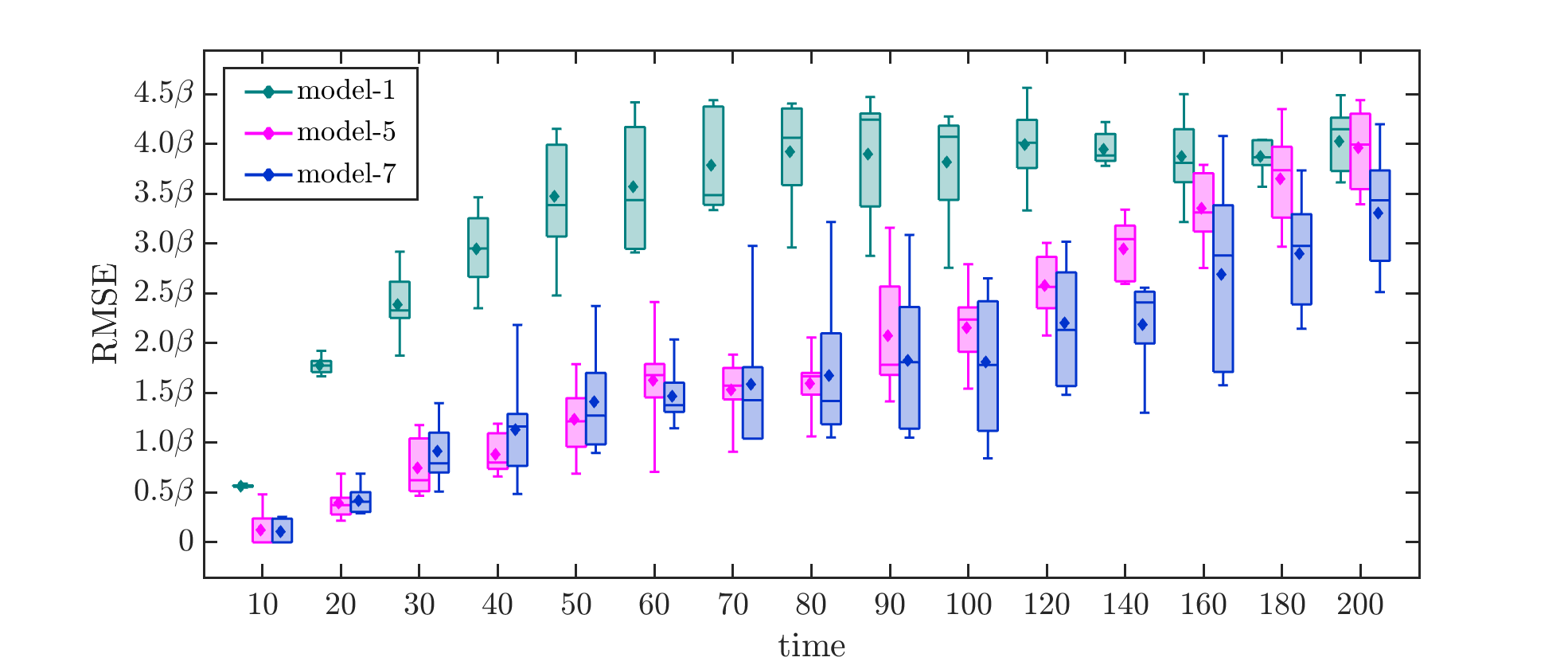}
    \caption{RMSEs achieved by model-1, model-5, and model-7 for different time periods of the planning horizon.}
    \label{fig:RMSE1}
\end{figure}

Looking at the inventory profiles shown in Figure \ref{fig:T200}, we notice that the accuracy of the model predictions seems to decrease in later time periods. We confirm this by plotting the RMSEs computed for smaller time intervals (10 time periods per time interval up to time 100 and 20 time periods per time interval after time 100) in Figure \ref{fig:RMSE1}, where the increasing trend in RMSE is particularly clear for model-5 and model-7. One way to interpret this trend is stated in the following hypothesis.

\begin{hypothesis}
    Planners focus more on the earlier time periods, i.e. care less about being optimal in their plan for the later time periods. \textup{-- True. Given that actual orders are only available for the first few time periods and the opportunity to re-optimize the production plan on a regular basis, it is not worth making a big effort in optimizing the later time periods where the actual demand will almost certainly be different from the current forecast.}
\end{hypothesis}

\begin{figure}[ht!]
\captionsetup{}
    \centering
    \includegraphics[width=15cm]{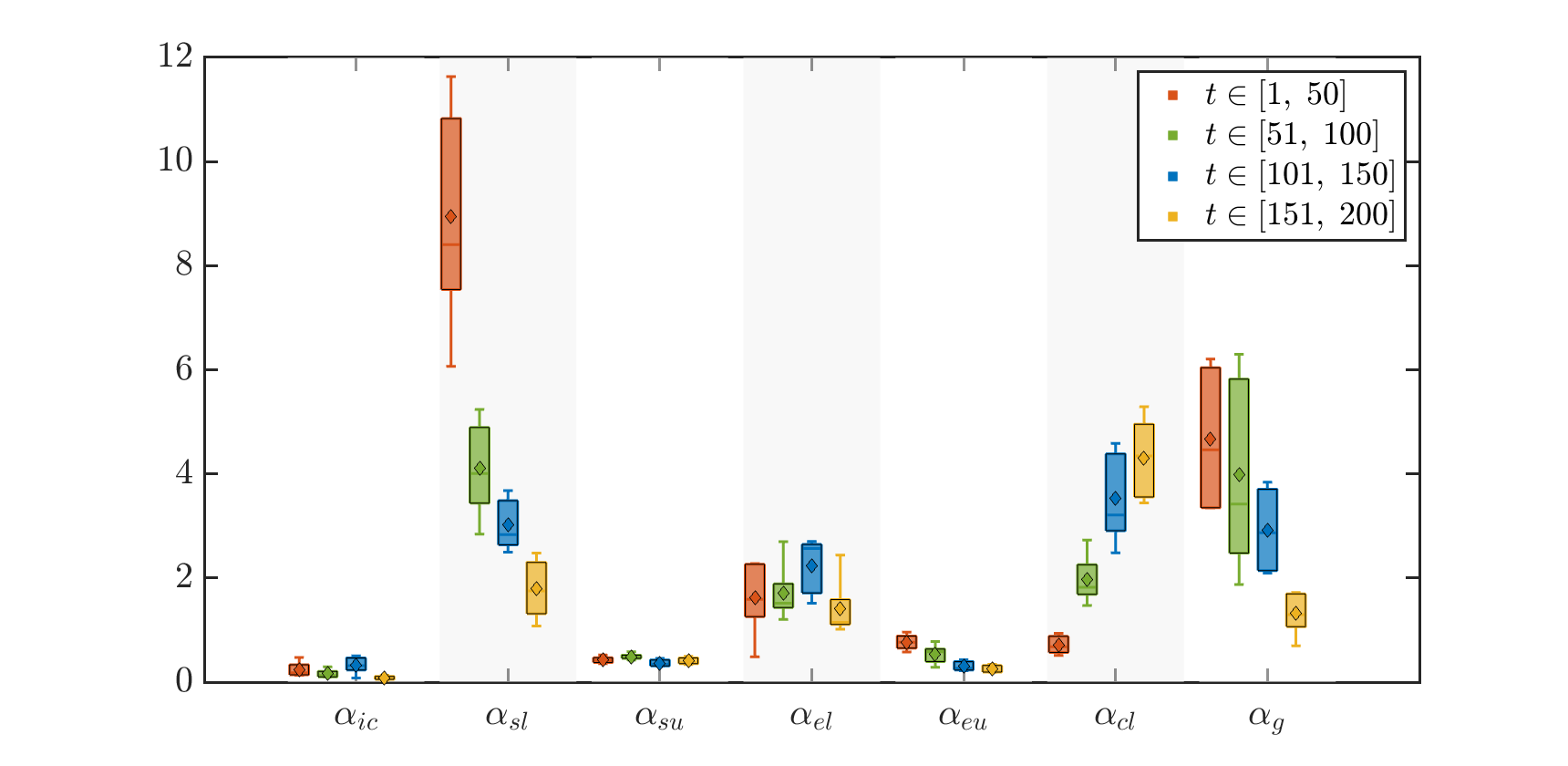}
    \caption{Estimates of the time-dependent weighting factors across the four 50-period time buckets obtained for model-7-t.}
    \label{fig:tim_costpref}
\end{figure}

Inspired by the time-varying prediction accuracy, we now want to explore whether planners also emphasize the various objectives differently in different time periods. To this end, we divide the planning horizon into four time buckets, each consisting of 50 time periods, and allow $\alpha$ to be different for each time bucket. The model trained using these time-dependent weighting factors is denoted by model-7-t. Figure \ref{fig:tim_costpref} shows the obtained estimates for $\alpha$ across the four time buckets.

First, we observe that $\alpha_{sl}$ decreases over time, which leads to the following hypothesis.

\begin{hypothesis}
    Planners try to avoid low inventory more in the earlier time periods. \textup{-- True. Actual demands are realized in the earlier time periods, and when they are higher than previously forecasted, backlogs can occur, which do not only affect the service level for the current orders but also the ability to fulfill orders in later time periods. In that case, a strong emphasis is placed on scheduling the production such that sufficiently high inventory levels can be maintained.}
\end{hypothesis}

We see an opposite trend for $\alpha_{cl}$, i.e. it increases over time.

\begin{hypothesis}
    Maintaining consistent cycles times is more important in the later time periods. \textup{-- True. In the earlier time periods, the main focus is on fulfilling the current orders, where deviating form the targeted cycle length is perfectly fine. When the demand is higher than forecasted, the cycle length is increased to increase production. Conversely, when the demand is lower than forecasted, the cycle length is reduced to avoid excessive inventory. Keeping the cycle lengths within the desired range is more important in later time periods as this is one main way to address demand uncertainty.}
\end{hypothesis}

Interesting, the results of model-7-t show that $\alpha_{sl}$ is greater than $\alpha_{el}$, especially in the earlier time buckets; this is the opposite of what we observed in model-7 (cf. Figure \ref{fig:costpref}).

\begin{hypothesis}
    Planners pay more attention to the inventory level at the start rather than at the end of a campaign. \textup{-- True. This is in line with the rationale provided to invalidate Hypothesis \ref{hyp:el_more_important}.}
\end{hypothesis}

\begin{figure}[ht!]
    \centering
    \subfloat{\includegraphics[width=5.0cm]{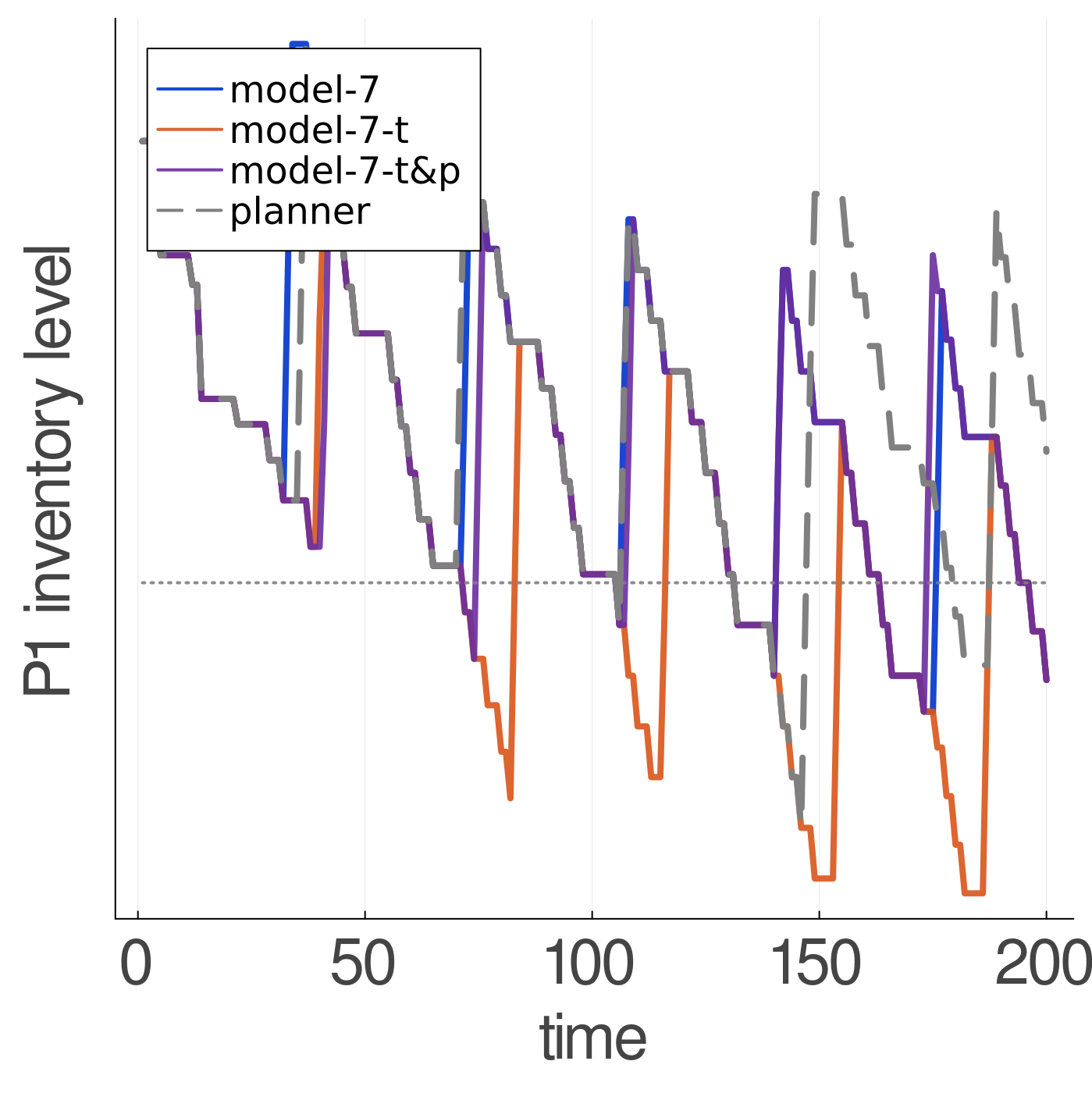}}
    \subfloat{\includegraphics[width=5.0cm]{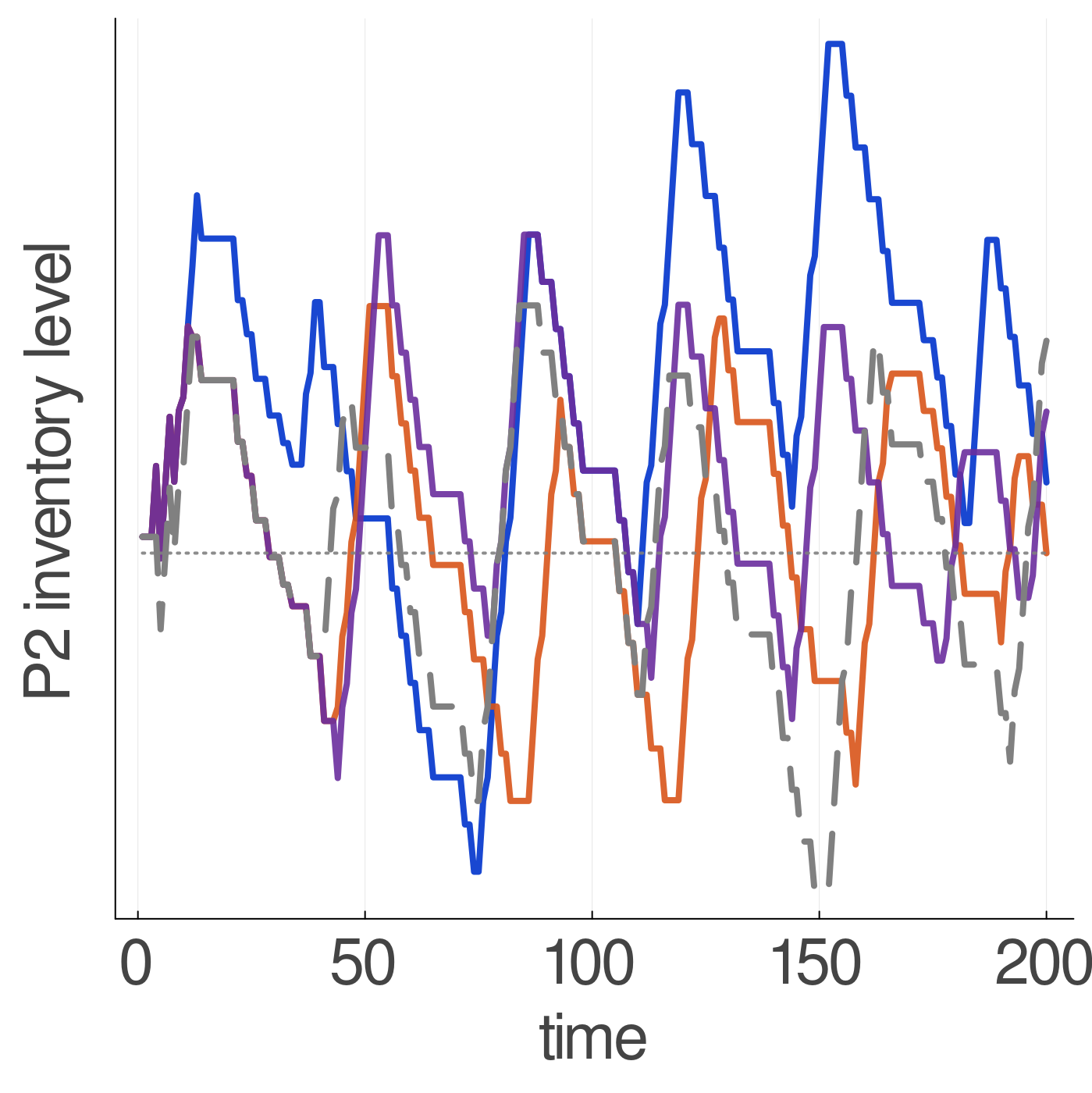}}
    \subfloat{\includegraphics[width=5.0cm]{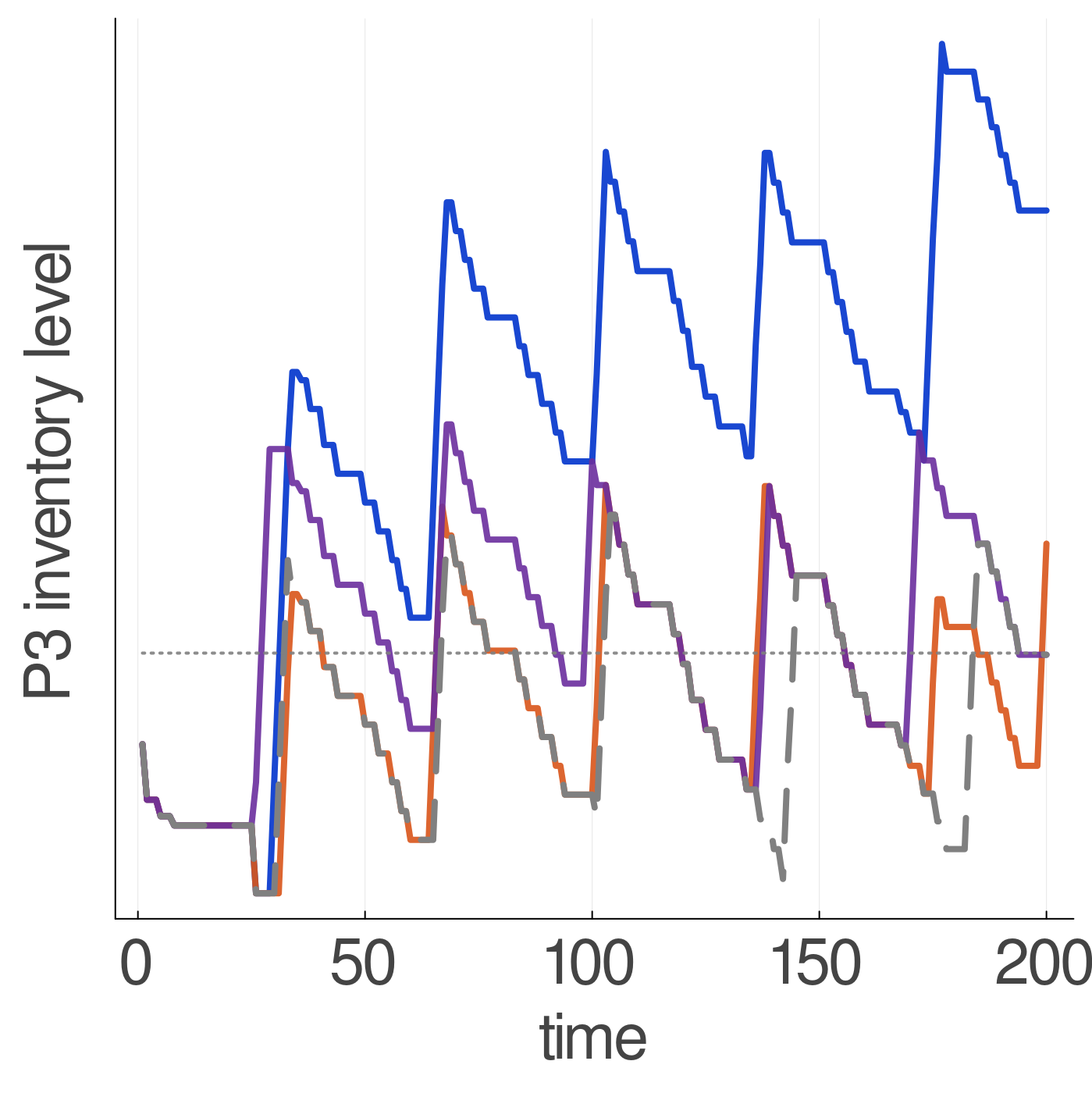}} \\
    \subfloat{\includegraphics[width=5.0cm]{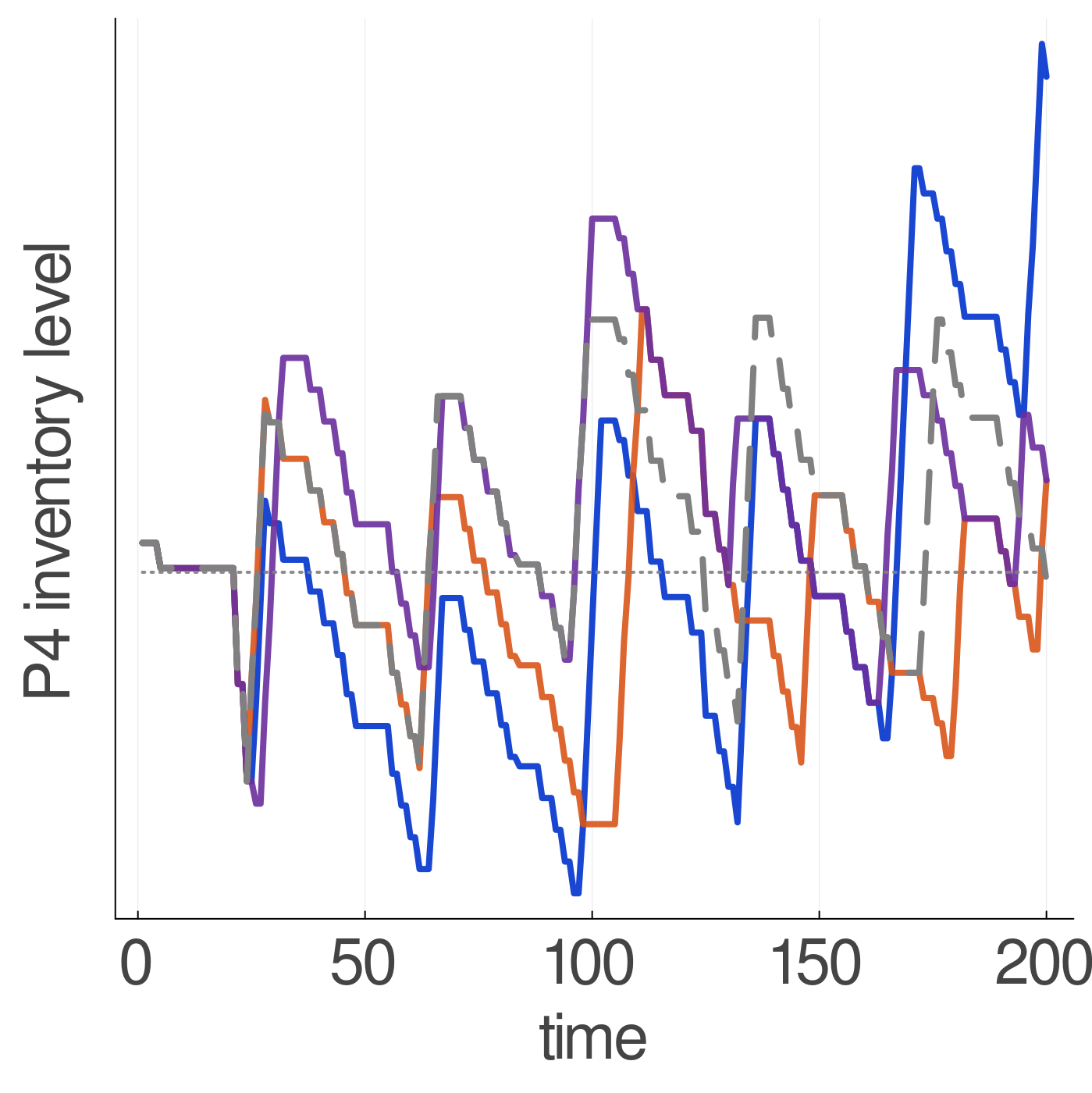}}
    \subfloat{\includegraphics[width=5.0cm]{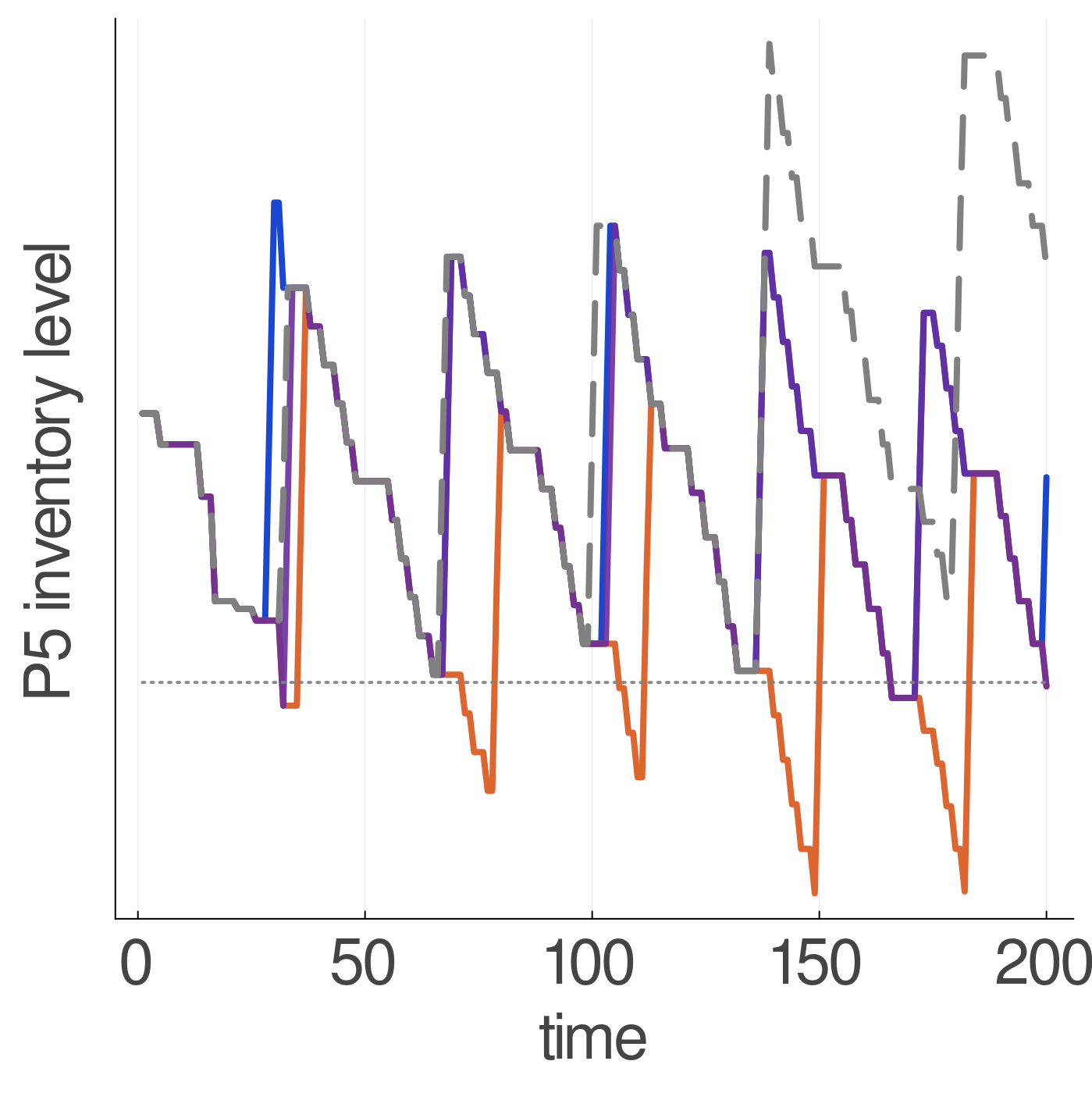}}
    \subfloat{\includegraphics[width=5.0cm]{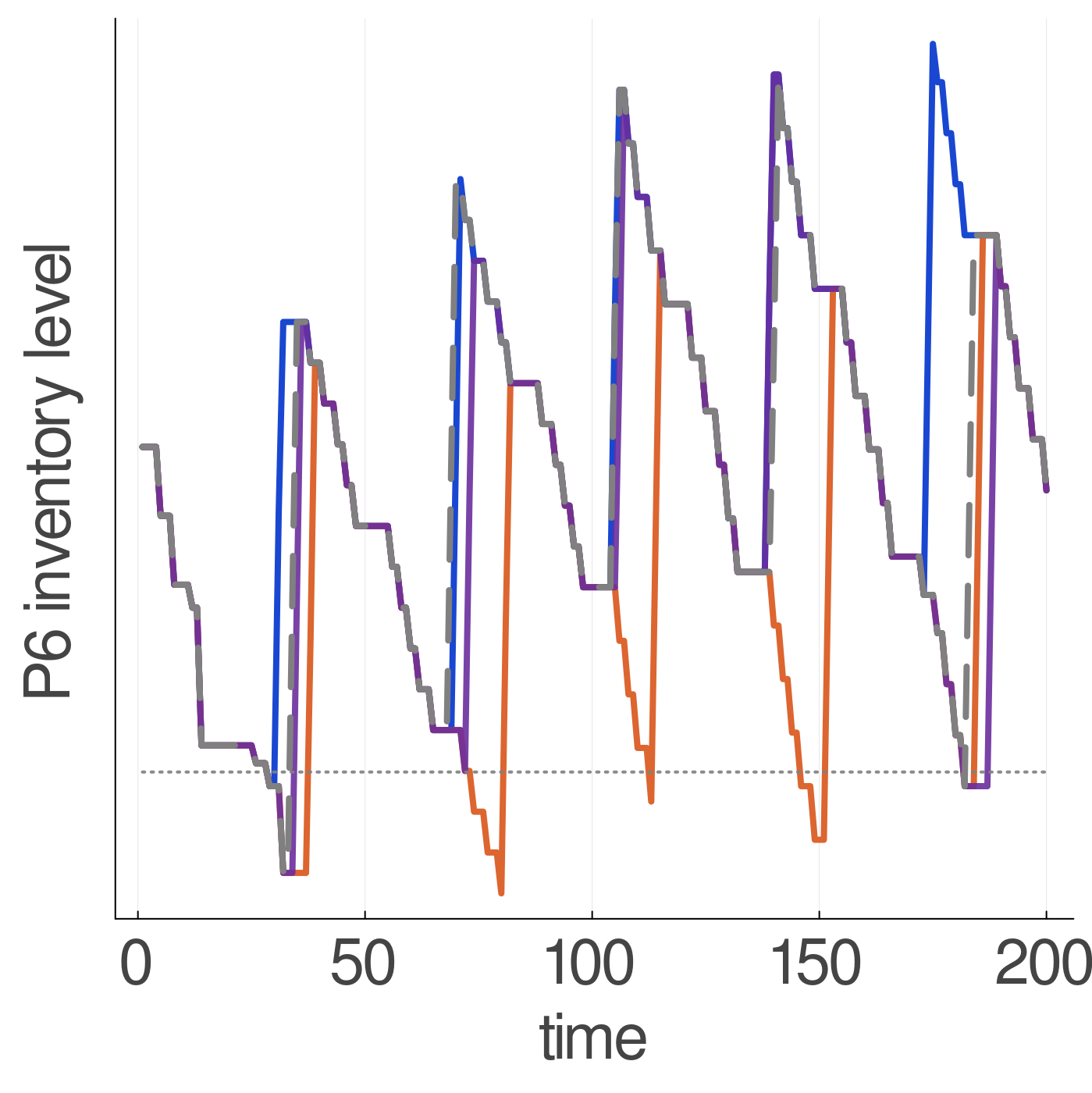}} \\
    \subfloat{\includegraphics[width=5.0cm]{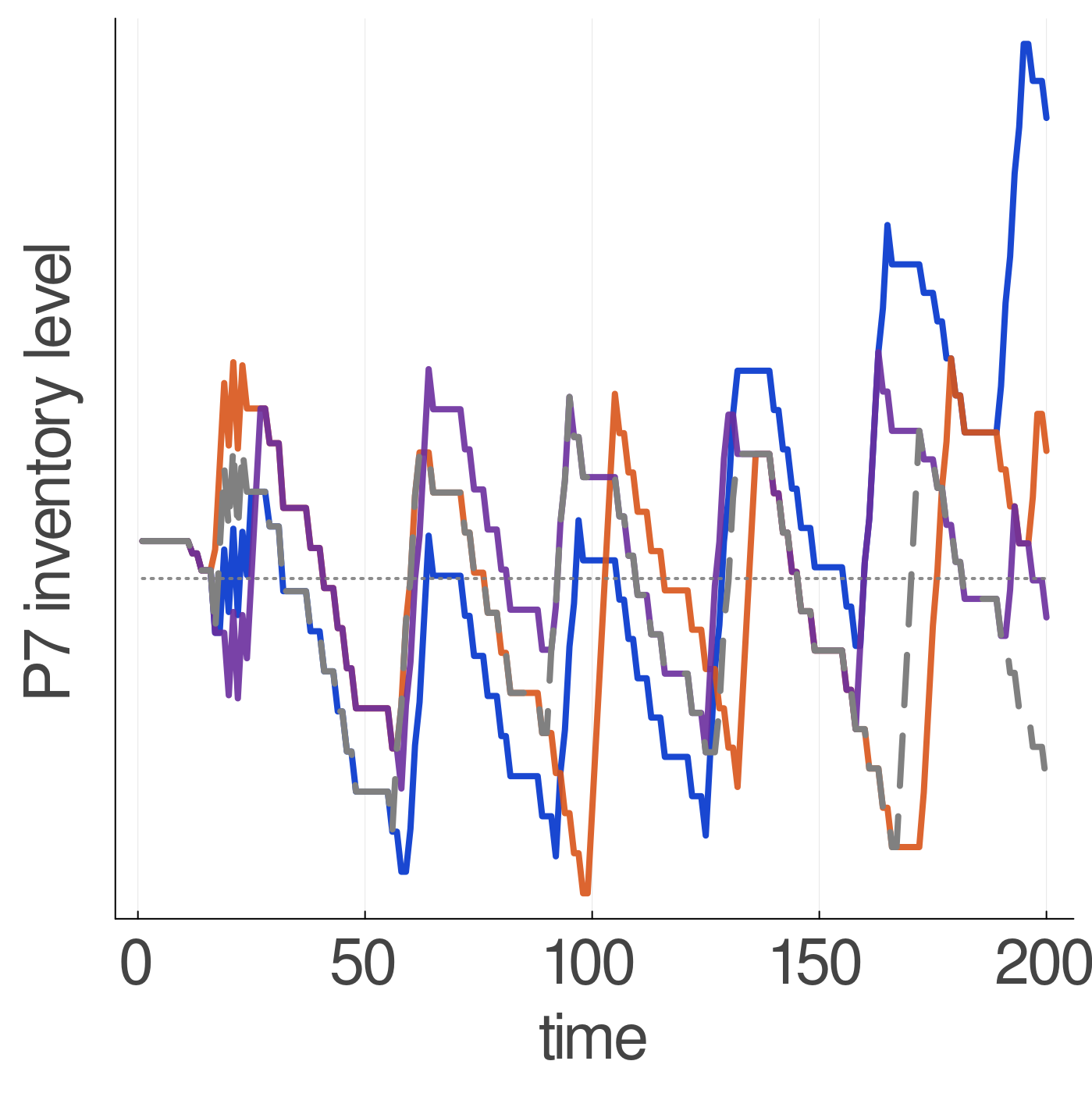}}
    \subfloat{\includegraphics[width=5.0cm]{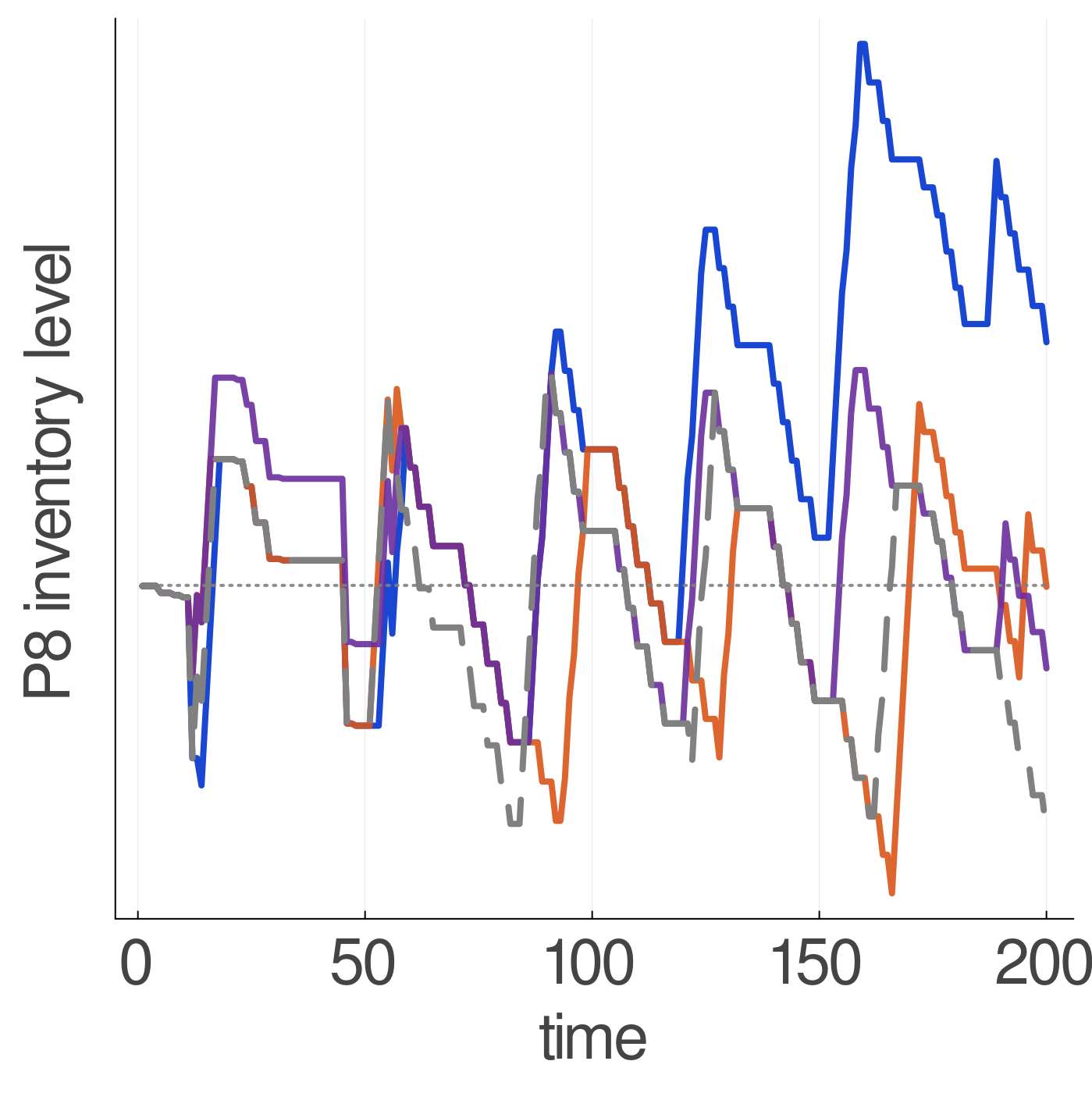}}
    \caption{Comparing for a particular instance the actual inventory profiles for all eight products with the profiles predicted by model-7, model-7-t, and model-7-t\&p.}
    \label{fig:T200_prodtim}
\end{figure}

\begin{figure}[ht!]
    \centering
    \includegraphics[width=15cm]{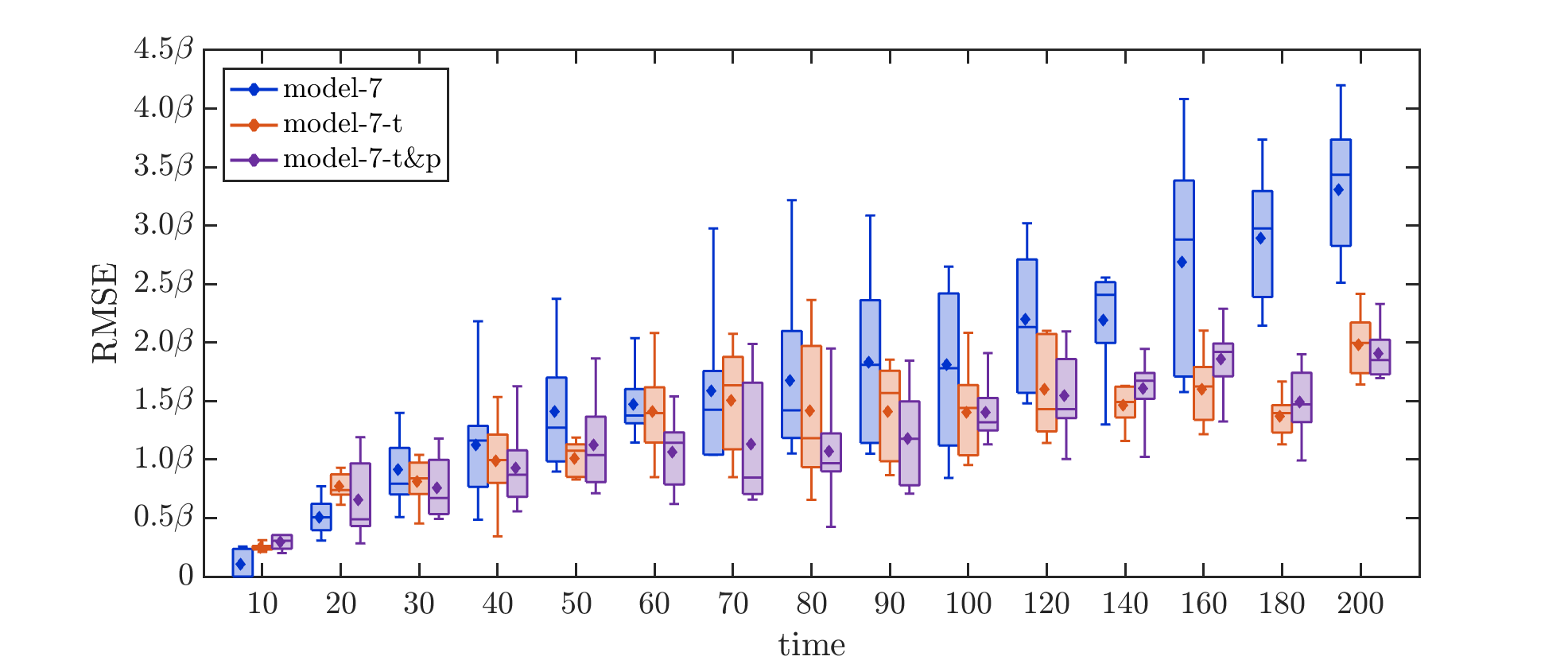}
    \caption{RMSEs achieved by model-7, model-7-t, and model-7-t\&p for different time periods of the planning horizon.}
     \label{fig:RMSE3}
\end{figure}

We also observe that the prediction accuracy of model-7 and model-7-t can vary significantly across different products (see Figure \ref{fig:T200_prodtim}). This has motivated us to build another model in which the objective weighting factors are also product-dependent, which we call model-7-t\&p. In the example shown in Figure \ref{fig:T200_prodtim}, we can see that model-7-t\&p predicts the inventory profiles for certain products, e.g. P5, more accurately than model-7-t. The is also true for the overall prediction accuracy, as shown in the RMSE plot in Figure \ref{fig:RMSE3}.

\begin{figure}[ht!]
    \centering
    \includegraphics[width=\textwidth]{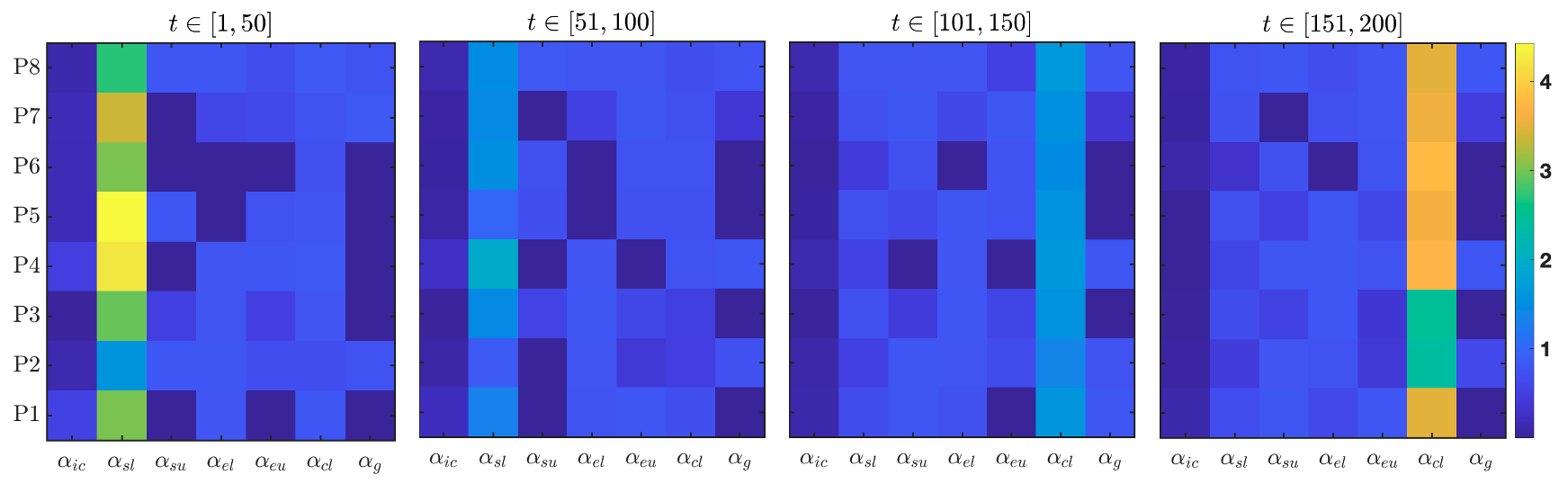}
    \caption{Heat maps showing the average values of $\alpha$ in model-7-t\&p for all four time buckets and eight products.}
    \label{fig:heatmap}
\end{figure}

The heat maps in Figure \ref{fig:heatmap} show the mean values of $\alpha$ for all time buckets and products. We again observe that $\alpha_{sl}$ and $\alpha_{cl}$ take large values toward the beginning and end of the planning horizon, respectively. Moreover, we now see that the weights can also be significantly different across the eight  products.

\begin{hypothesis}
\label{hyp:product_dependence}
    Planners tend to focus on some products more than others. \textup{-- True. Different products can have different profit margins, and different customers can also have different priorities. Planners typically try to maintain a higher inventory level for important products to ensure on-time fulfillment of high-priority orders.}
\end{hypothesis}

Regarding Hypothesis \ref{hyp:product_dependence}, we should mention that it has not been revealed to us whether the product-dependent weights we obtained align with the true relative importance of the different products or not. We also do not want to make such specific conclusions from these results since there can be many other factors that lead to these product-dependent weights. For example, a large $\alpha_{sl}$-value may also be indicative of a significant and consistent underforecasting of demand, which increases the likelihood of understock and results in an even greater need to maintain a sufficiently high inventory level. Although no clearer insights can be drawn at this point, we plan in our future work to use our results as a basis for further discussions with the planners to untangle the various potential influencing factors and further refine our model.

\subsection{Computational analysis}

\begin{figure}[!htbp]
     \centering
     \includegraphics[width=12cm]{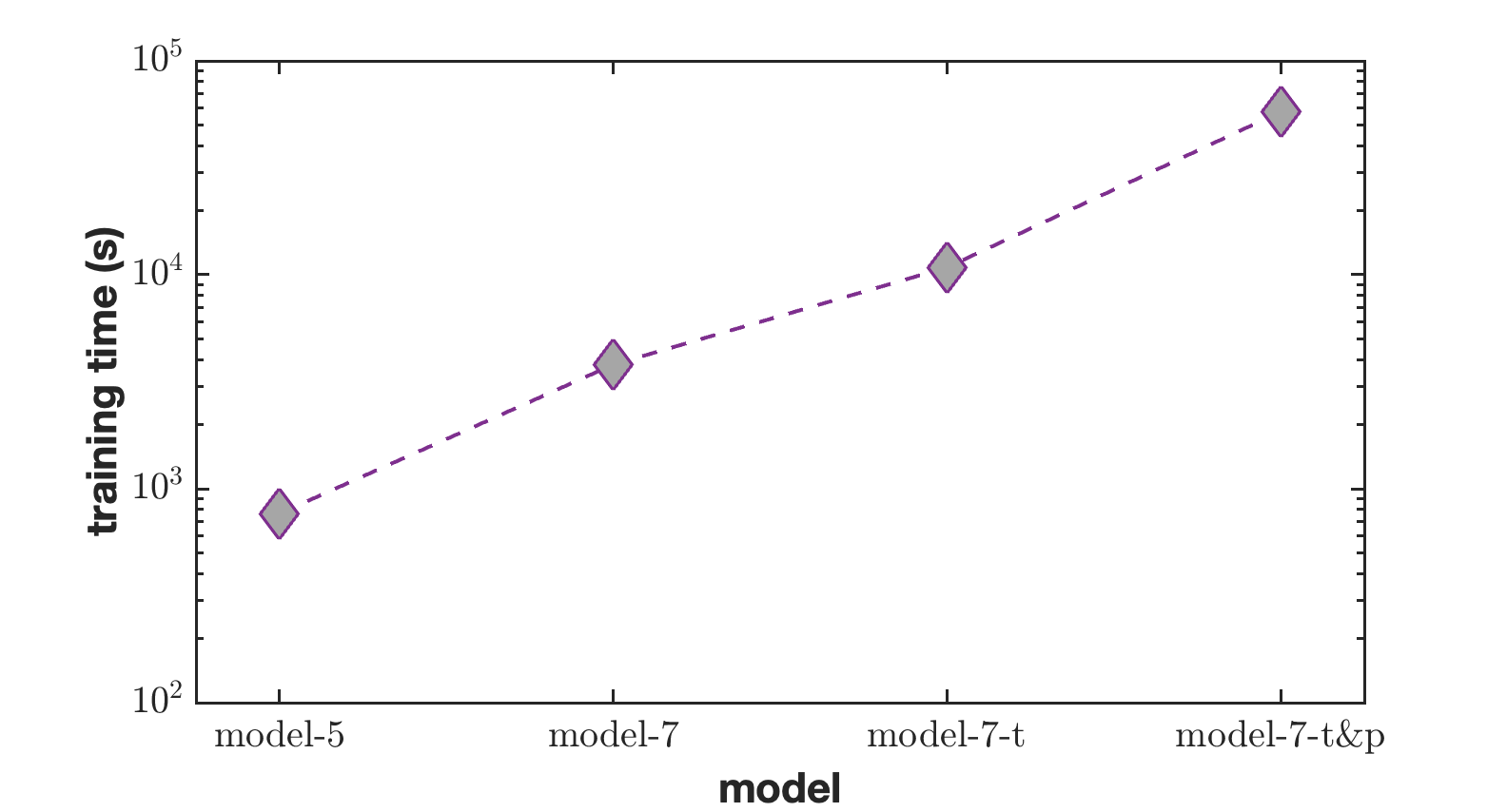}
     \caption{Training times required by the proposed cutting-plane algorithm for different models.}
      \label{fig:train_time}
\end{figure}

We conclude our analysis by examining the computational performance of the proposed framework across all four models in terms of training time and convergence behavior of the cutting-plane algorithm. Figure \ref{fig:train_time} shows the training time for each model on a logarithmic scale. Model-5, which estimates five objective weighting factors based solely on inventory-related penalties, requires approximately 800 seconds to train. Introducing the two additional parameters for cycle length deviations and production gaps in model-7 increases the training time by about one order of magnitude to around 3000 seconds. Extending the model to consider time-dependent weighting factors in model-7-t further raises the training time to approximately 10,000 seconds, as the number of parameters grows proportionally with the number of time windows. Finally, model-7-t\&p, which additionally allows for product-dependent weights, incurs the highest training cost at nearly 60,000 seconds, reflecting the combinatorial growth in parameters across both time windows and products. Despite this increase, it is important to note that all models are trained offline; once the cost coefficients are estimated, production planning decisions can be generated by simply solving problem (\ref{eqn:optproblem}) without any additional training overhead, making the approach practical even for the most complex model variant.

\begin{figure}[!htbp]
     \centering
     \includegraphics[width=10cm]{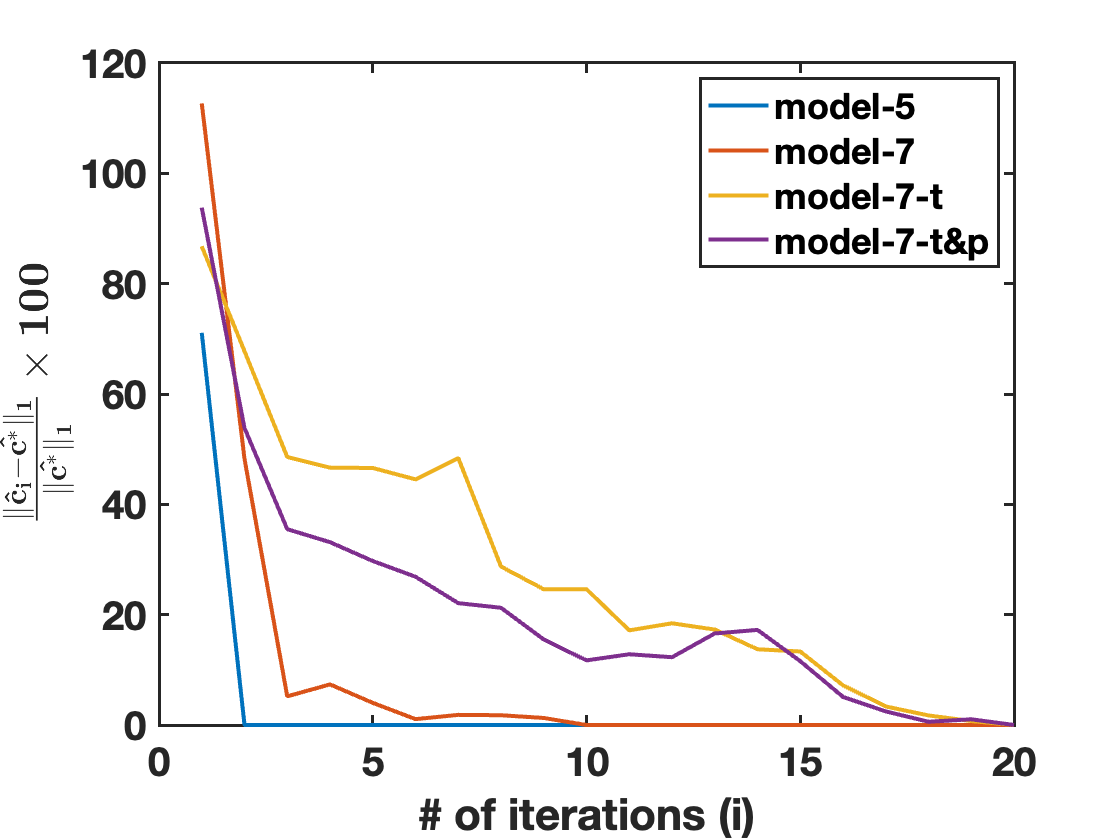}
     \caption{Convergence of the IO algorithm for different models.}
      \label{fig:convergence}
\end{figure}

Figure \ref{fig:convergence} illustrates the convergence behavior of the cutting-plane algorithm for solving problem \eqref{eqn:IOP-SL} across all four models, where the y-axis reports the normalized cost coefficient error  as a function of the iterations. Here, model-5 converges most rapidly, reaching near-zero error within just a few iterations, consistent with its relatively small number of parameters and the limited complexity of the corresponding cut-generating subproblems. Model-7 follows a similar trajectory, converging within approximately ten iterations. In contrast, model-7-t and model-7-t\&p exhibit slower convergence, with the normalized error reaching near-zero by iteration 20. This behavior is expected since the larger number of parameters increases the computational complexity of both the master problem and the cut-generating subproblems. Nonetheless, all models converge within the allotted 20 iterations, demonstrating the practical efficacy of the cutting-plane approach.

\section{Conclusions} 
\label{sec:Conclusions}

In this work, we were given a production planning optimization problem for which the constraints were known but the objective function was not. Assuming that the right objective function is captured in the decisions of expert planners that have been solving this production planning problem for many years, we applied a data-driven IO approach to learn the unknown objectives of these planners from production plans that they have generated in the past. Specifically, we formulated the production planning problem as an MILP with a weighted-sum objective function that considers multiple objectives which planners could be considering in their decision-making; the resulting IO problem was to infer the corresponding weighting factors. 

In the industrial case study provided by Dow, we were able to achieve reasonable prediction accuracy with only 50 training data points, each representing a historical production plan. This level of data efficiency can be attributed to the fact that all production planning constraints were explicitly captured in the IO framework, which means that inference could focus on the part that is truly unknown, namely the cost coefficients. Data efficiency is indeed important from a practical standpoint since this kind of data is not easy to come by in this particular production planning context.

Perhaps even more importantly, the results also highlight the interpretability of the proposed approach. Since the inferred weighting factors are directly associated with the corresponding objectives, they provide immediate insights into the relative significance between these objectives as perceived by the expert planners. This allowed us to formulate specific hypotheses about the planners' decision-making process and discuss them with an expert planner, whose feedback not only helped confirm or invalidate our hypotheses but also enabled further refinement of the model. Ultimately, we obtained a model that not only predicts the planners' decisions well but is also interpretable such that planners can see how it aligns with the rationale in their own decision-making. We believe that only such a model can serve as a decision-support tool that is trusted by the user.

\section*{Acknowledgments}

The authors gratefully acknowledge the financial support from the National Science Foundation under Grant \#2044077 as well as the Minnesota Supercomputing Institute (MSI) at the University of Minnesota for providing resources that contributed to the research results reported in this paper. We would also like to thank Chad Dandron, Ellen Du, Kyle Harshbarger, Ellen Murray, and Jeff Tazelaar for their support and advice on the project.

\bibliographystyle{abbrvnat}
\bibliography{example_paper, library, library_rg}

\end{document}